\documentclass[review]{elsarticle}

\usepackage{lineno,hyperref}
\modulolinenumbers[5]
\usepackage{amsfonts,amssymb,amsmath,amstext,multirow}
\usepackage[english,activeacute]{babel}
\usepackage{graphics,graphicx}
\usepackage{soul} 
\usepackage{cancel} 
\usepackage{color}

\newtheorem{theorem}{Theorem}
\newtheorem{corollary}[theorem]{Corollary}
\newtheorem{lemma}{Lemma}

\newtheorem{remark}{Observation}

\usepackage{listings}

\begin{document}
\begin{frontmatter}


\title{ A Petrov-Galerkin multilayer discretization to second order elliptic boundary value problems}


\author[add1]{Tom‡s Chac\'{o}n Rebollo}
\ead{chacon@us.es}
\author[add2]{Daniel Franco Coronil\corref{cor1}}
\ead{franco@us.es}
\author[add3]{Fr\'ederic Hecht}
\ead{hecht@ann.jussieu.fr}
\cortext[cor1]{Departamento EDAN, Universidad de Sevilla, C/Tarfia, s/n, 41012 Sevilla (Spain)}
\address[add1]{Departamento EDAN \& IMUS, Universidad de Sevilla, C/Tarfia, s/n, 41012 Sevilla (Spain)}
\address[add2]{Departamento EDAN, Universidad de Sevilla, C/Tarfia, s/n, 41012 Sevilla (Spain)}
\address[add3]{Laboratoire Jacques-Louis Lions, Universit\' e de Paris VI. 5 Place Jussieu, 75005 Paris (France)}


\begin{abstract}
We study in this paper a multilayer discretization of second order elliptic problems, aimed at providing reliable multilayer discretizations of shallow fluid flow problems with diffusive effects. This discretization is based upon the formulation by transposition of the equations. It is a Petrov-Galerkin discretization in which the trial functions are piecewise constant per horizontal layers, while the trial functions are continuous piecewise linear, on a vertically shifted grid.

We prove the well posedness and optimal error order estimates for this discretization in natural norms, based upon specific inf-sup conditions.

We present some numerical tests with parallel computing of the solution based upon the multilayer structure of the discretization, for academic problems with smooth solutions, with results in full agreement with the theory developed.
\end{abstract}

\end{frontmatter}

\section{Introduction and motivation}
This paper deals with the numerical approximation of the Poisson and related problems by means of layer-wise discontinuous solutions. It is motivated by the construction of multilayer discretizations of fluid flow equations on sha\-llow domains.  This produces a reduction of the dimensionality of the approximated problem, yielding in practice a domain decomposition discretization that is solved by parallel procedures. \par
This is the case for instance of the works of Fern‡ndez-Nieto and co-workers \cite{bonaventura:2018, fernandezNieto:2016} and of Sainte-Marie and co-workers \cite{SainteMarie} to solve the Navier-Stokes equations.  In these papers the velocity is discretized by layer-wise smooth functions, discontinuous across the layer boundaries. The conservation of momentum needs the continuity of diffusive fluxes across layer boundaries. This is ensured in these papers by specific techniques based upon finite-difference discretizations in the vertical direction. Both papers perform stability analysis for the discretizations introduced, however no convergence proofs are reported. This mainly occurs because the discretizations are purely finite-difference approximations, that are not linked to variational formulation of the targeted equations. 
\par In this paper we look for stable multilayer discretizations for which convergence proofs are reachable. More concretely, we look for stable discretizations that are related to variational formulation of the equations. A possible discretization meeting these criteria could be provided by the Discontinuous Galerkin (DG) method. This method was introduced in the 90's by Cockburn and Shu for hyperbolic conservation laws (cf. \cite{cocshu}), and later adapted to elliptic problems, by re-writing as two coupled first-order conservation laws for the unknown and its gradient (cf. \cite{cockburn, pietro} and references therein).  DG methods for elliptic problems are based upon the flux formulation of the elliptic problem, that includes both the original unknown so as its gradient as unknowns of the discretization. This allows to approximate both the unknown and its gradient by discontinuous finite elements. The stability of the formulation strongly depends on the choice of the numerical fluxes of both unknowns.
\par An alternative discretization could be provided by the Discontinuous Petrov Galerkin  methodology, introduced in a series of papers by Sacco and co-workers (cf. e.g. \cite{botasso, causin}) and systematically studied by Demkowicz and Gopalakrishnan (see \cite{demko1, demko2, demko3, demko} and references therein). The DPG methodology approximates the unknown by either continuous or discontinuous trial functions, while the tests functions necessarily need to be discontinuous. This formulation admits three equivalent interpretations: a Petrov-Galerkin method with test functions that achieve the supremum in the inf-sup condition, a minimum-residual method with residual measured in a suitable dual norm, and a mixed formulation where one solves simultaneously for the Riesz representation of the residual. The stability of the formulation is a direct consequence of the first of these three interpretations.
\par Here we propose a discretization specifically adapted to multilayer discretizations, for elliptic problems related to the Poisson equations. The main idea is to start from the solution by transposition of the equations, much as the ultra-weak formulation considered in the works \cite{botasso, causin} by Sacco and co-workers. The solution by transposition naturally belongs to $L^2$ spaces, and thus admits piecewise discontinuous approximations. Based upon this procedure, we propose a Petrov-Galerkin discretization for cylindrical domains, in which the solution is a layer-wise constant function, while the test functions are continuous piecewise affine polynomials in the vertical direction, with knots in a shifted grid. We derive a single \lq\lq recepy" to build this kind of discretizations: Approximate the vertical derivative of the unknown by the vertical derivative of its interpolate on the test functions space. 
\par We prove the well-posedness of this discretization, given by an inf-sup condition satisfied by the bilinear form appearing in the discrete problem. Also, we prove optimal order error estimates for smooth solutions. We further extend the multilayer discretization technique introduced to domains which are vertical deformations of cylindrical domains, which appear in shallow water problems when a flat surface is deformed. We also prove well-posedness and optimal order error estimates for this extension.  We further extend the method to Neumann boundary conditions, by slightly changing the test space. 
\par
We finally present some 3D numerical tests by parallel solution of the resulting linear system, taking advantage of the multilayer structure of the discretization, what allows to only solve 2D linear systems. These tests present speeds-ups rates ranging from 20 to 50, while presenting optimal convergence orders for smooth functions, in agreement with the theoretical expectations.
\par 
The paper is organized as follows. Section \ref{se:multilayer} introduces the motivation and the basic multilayer discretization that we consider, for the Poisson problem in a cylindric domain. This discretization is studied in Section \ref{wellpodness}, where stability, convergence and obtention of error estimates are analyzed.  In Section \ref{se:nonflat} the discretization is extended to Poisson problems in domains with non-flat upper boundary, also analyzing stability, convergence and error estimates. Section \ref{se:neumann} is devoted to the extension of the discretization to Neumann boundary conditions. Finally Section \ref{se:tests} present several 3D numerical tests for each of the cases considered, for smooth solutions.
  
\section{Multilayer approximation}\label{se:multilayer}
Let us consider a cylindrical domain $\Omega \,=\,\omega\:\times \:(0,L)$ where $\omega\,\subset\,\mathbb{R}^d$ is bounded domain, and $d\ge 1$ is an integer number. Let us consider the homogeneous Dirichlet Poisson problem as a model problem:  \\ 
\begin{equation}\label{eq:poispb0}
\left\{\begin{array} {ccccc}
-\Delta v &= & f & \mbox{in} & \Omega, \\ [2ex]
v & = & 0 & \mbox{on} & \partial \Omega.
\end{array}\right.
\end{equation}
As is standard, this problem can be written in variational form: Given $f \in H^{-1}(\Omega)$, find $v \in H_0^1(\Omega)$  such that
\begin{equation}\label{fv:poispb}
a\:(v,\varphi) \,=  \,L\:(\varphi),\qquad \forall \;\varphi \,\in\,H_0^1(\Omega), 
\end{equation}
with 
\begin{equation} \label{forma}
a\:(v,\varphi) \,=\,\displaystyle\int_\Omega \nabla v \:\nabla \varphi\, d\Omega, \,\mbox{ and   }L\:(\varphi)\,=\,<f , \varphi>_\Omega.
\end{equation}
 
The solution of this problem is also the solution of the transposition formulation of problem \eqref{eq:poispb0}, given by  (Cf. \cite{brezis}): Find $v \in L^2(\Omega)$  such that for all $\phi \;\in\; L^2(\Omega)$,

\begin{equation}\label{eq:poistrans2}
\int_\Omega v \phi \,d\Omega\,=\,L\:({\cal T} \phi),
\end{equation}
where ${\cal T} \phi =\varphi$  is the solution of
\begin{equation}\label{eq:poistransaux}
\left\{\begin{array} {ccccc}
-\Delta \varphi &= & \phi & \mbox{in} & \Omega, \\ [1ex]
\varphi & = & 0 & \mbox{on} & \partial \Omega.
\end{array}\right.
\end{equation}

Problem (\ref{eq:poistrans2}) admits a unique solution in $L^2(\Omega)$, that necessarily coincides with the solution of \eqref{fv:poispb}. Also, problem (\ref{eq:poistrans2}) admits a unique solution when $f$ is less smooth (actually, when $f \in (H^1_0(\Omega)\cap H^2(\Omega))'$) whenever the problem (\ref{eq:poistransaux}) is regular, in the sense that $\varphi \;\in\; H^1_0(\Omega)\cap H^2(\Omega)$ (Cf.~\cite{brezis}).
\par
\begin{figure}
  \begin{center}
   \makebox[\linewidth]{\includegraphics[height=.6\linewidth]{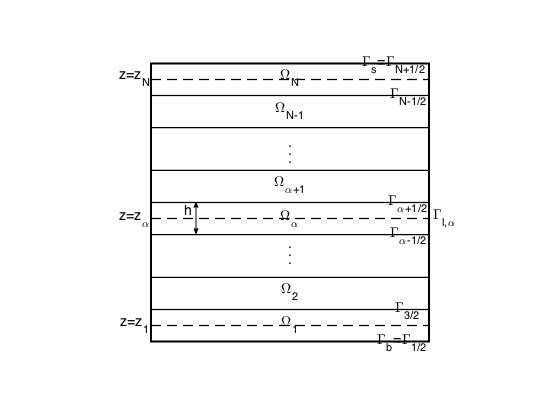}}
  \caption{Splitting of the domain $\Omega$ \label{fig:dom}}
  \end{center}
\end{figure}

 However, we are not interested in approximating the formulation by transposition \eqref{eq:poistrans2}, rather we use it as a base to build our multilayer discretization of the weak standard formulation \eqref{fv:poispb}. To do it, let us split the domain $\Omega$ along the vertical direction into $N\,\geq\,1$  plane layers of constant thickness $h=L/N$ with $N + 1$ interfaces $\Gamma_{\alpha+\frac{1}{2}}$ of equations $z \,=\, z_{\alpha+\frac{1}{2}}$ for $\alpha = 0$, 1, . . . , $N$ (see  Figure~\ref{fig:dom}). Hence, we split $\overline\Omega=\displaystyle\bigcup_{\alpha=1}^N \overline{\Omega_\alpha}$ with
$$
\Omega_\alpha =\left\{(x,z);\;x\,\in\,w\;\mbox{and}\; z_{\alpha-\frac{1}{2}}\,<\,z\,<\,z_{\alpha+\frac{1}{2}}\,\right\}.
$$

We assume  $\partial\Omega = \Gamma_b \cup \Gamma_s \cup \Gamma_l$, where $\Gamma_b$ and $\Gamma_s$ are the boundaries of the bottom and the top, respectively, and $ \Gamma_l$ is the vertical boundary of the domain. We shall denote by $z_b = z_{\frac{1}{2}}=0$ and $z_s = z_{N+\frac{1}{2}}=L$ the equations of the bottom and the top interfaces $\Gamma_b$ and $\Gamma_s$, respectively.  Notice that  $\partial\Omega_\alpha = \Gamma_{\alpha-\frac{1}{2}} \cup \Gamma_{\alpha+\frac{1}{2}} \cup\Gamma_{l,\alpha}$, with
$$
\Gamma_{l,\alpha} =\left\{(x,z)\in \mathbb{R}^d;\;x\,\in\,\partial w\;\mbox{and}\; z_{\alpha-\frac{1}{2}}\,\leq\,z\,\leq\,z_{\alpha+\frac{1}{2}}\,\right\}.
$$

Then the vertical boundary is split as $\Gamma_{l}=\displaystyle\bigcup_{\alpha=1}^N {\Gamma_{l,\alpha}}$. 

Also, we denote $z_\alpha=\displaystyle\frac{z_{\alpha-\frac{1}{2}}\,+\,z_{\alpha+\frac{1}{2}}}{2}$, $\forall\;\alpha = 1$, 2, . . . , $N$.

We shall assume from now on that $\omega$ is polygonal. Let us now introduce discretization space for the unknown $v$,
$$X_h= \displaystyle\left\{v_h\,\in\,L^2(\Omega)\, /\; v_h=\sum_{\alpha=1}^N v_h^\alpha \otimes \kappa_\alpha\, \mbox{with  } v_h^\alpha\in\,V_{0,k}(\omega)\right\},
$$  
where $\kappa_\alpha$ is the characteristic function of the interval $(z_{\alpha -1/2},z_{\alpha+1/2})$ and $V_{0,k}(\omega)$ is a finite element sub-space of $H^1_0(\omega)$. The sub-index $k$ denotes the horizontal grid size. We assume that $V_{0,k}(\omega)=V_k(\omega)\,\cap\,H_0^1\:(\omega)$ where $V_k(\omega)$ is a Lagrange finite element space, constructed when $d=2$ with either triangular or quadrilateral elements, i. e., there exists a grid ${\cal M}_k$ such that
$$
V_k(\omega)=\{ w_h \in C^0(\overline{\omega})\, |\,w_h{|_K} \in R_l(K), \,\,\forall K \in {\cal M}_k\,\},
$$
where $R_l(K)= \left \{ \begin{array}{l} P_l(K) \mbox{ if }  {\cal M}_k \mbox{ is formed by triangles} \\
Q_l(K) \mbox{ if }  {\cal M}_k \mbox{ is formed by quadrilaterals},
\end{array}
\right .
$ for some integer $l \ge 1$; where $P_l(K)$ is the space of polynomials on $K$ of global degree less than or equal to $l$, and $Q_l(K)$ is the space of polynomials on $K$ of degree less than or equal to $l$ in each variable. When $d=1$, the triangulation is formed by segments $K$, and $R_l(K)=P_l(K)$.
\par

\begin{figure}
  \begin{center}
\makebox[\linewidth]{\includegraphics[height=.24\linewidth]{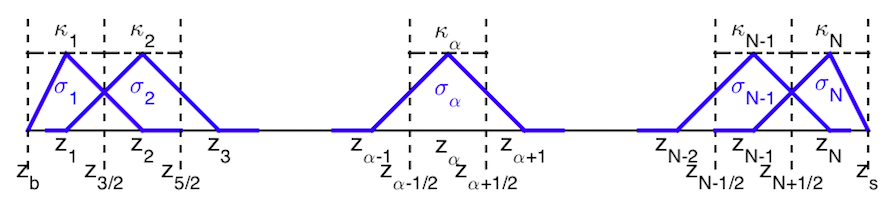}}
  \caption{Basis function on $X_h$ and $Y_h$ \label{fig:basis}}
  \end{center}
\end{figure}

We approximate the solution $v$ of problem \eqref{eq:poistrans2} by some $v_h \in X_h$. Observe that for all $\varphi\, \in \,H^2(\Omega)\cap H_0^1(\Omega)$,
\begin{equation}\label{eq:poistransalpha0}
\begin{array} {l}
\displaystyle-\int_\Omega v_h \:\Delta \varphi\, d\Omega\,=\,-\sum_{\alpha=1}^N\int_{\Omega_\alpha} v_h^\alpha \:\Delta \varphi\, d\Omega\,=\,
-\displaystyle\sum_{\alpha=1}^N\int_{\partial\Omega_\alpha} v_h^\alpha \partial_n\varphi d\Gamma\\ [2ex]+\displaystyle\sum_{\alpha=1}^N\int_{\Omega_\alpha} \nabla v_h^\alpha \nabla  \varphi d\Omega.
\end{array}
\end{equation}
As
$v_h^\alpha|_{\partial\omega}=0$ , then it holds
$$\int_{\partial\Omega_\alpha} v_h^\alpha \partial_n\varphi\,=\,\int_{\Gamma_{\alpha+\frac{1}{2}}} v_h^\alpha \partial_z\varphi\,-\,\int_{\Gamma_{\alpha-\frac{1}{2}}} v_h^\alpha \partial_z\varphi,$$
where $\partial_z$ denotes the vertical derivative. Thus,
$$
\displaystyle-\int_\Omega v_h \:\Delta \varphi\, d\Omega\,=\,a_h\:(v_h,\varphi), $$
where $a_h$ is a the bilinear form defined by
\begin{equation} \label{eq:formah}
\begin{array}{rcl}
a_h\:(v_h,\varphi) &=&\,\displaystyle\sum_{\alpha=1}^N (\nabla v_h^\alpha,\nabla  \varphi)_{\Omega_\alpha} + \displaystyle\sum_{\alpha=1}^{N-1}\:\int_{\Gamma_{\alpha+\frac{1}{2}}} (v_h^{\alpha+1}-v_h^\alpha) \partial _z\varphi\:  d\Gamma  \\  &&- \displaystyle\int_{\Gamma_s} v_h^N \partial _z\varphi \: d\Gamma + \displaystyle\int_{\Gamma_b} v_h^1 \partial _z\varphi\,  d\Gamma ;
\end{array}
\end{equation}
Next, we introduce the discrete test functions space, 
$$Y_h\,=\,\displaystyle\left\{\varphi_h\,\in\,L^2(\Omega)\, /\; \varphi_h\,=\,\sum_{\alpha=1}^N\:\varphi^\alpha\otimes\sigma_\alpha,\;\mbox{with}\;\varphi^\alpha\:\in\:V_{0,k}(\omega)\right\},$$
where $\sigma_1,\,\sigma_2,\,\cdots,\,\sigma_N$ are piecewise affine 1D functions on the intervals $[z_b, z_1]$, $[z_1,z_2], \cdots, [z_{N-1},z_N]$, $[z_N,z_s]$ that satisfy $\sigma_i(z_j)\,=\delta_{ij}$ and $\sigma_i(z_b)\,= \sigma_i(z_s)=0$ , for all $i, j = 1, 2, . . . , N$, as we see in Figure \ref{fig:basis}.
Then $Y_h$ is a subspace of $H_0^1(\Omega)$ formed by piecewise $\mathbb{P}_1$ polynomials in the variable $z$ with coefficients in $H^1_0(\omega)$.  Note that the intervals $[z_b, z_1]$, and $[z_N,z_s]$ have length $h/2$. Introducing them allows to have $dim (Y_h)=dim(X_h)$.

We now are in a position to state the multilayer discretization of (\ref{fv:poispb}) that we study in this paper,
\begin{equation}\label{eq:poistransalphadis}
\left\{\begin{array}{l}
\mbox{Find}\quad v_h\,\in\, X_h, \quad \mbox{such that}
\\[2ex]
a_h\:(v_h,\varphi_h) \,=\, L\:(\varphi_h), \qquad \forall \;\varphi_h\, \in \,Y_h.
\end{array}\right. 
\end{equation}
Notice that although the space $Y_h$ is not a sub-set of $H_0^1(\Omega) \cap H^2(\Omega)$, the action $a_h(v_h,\varphi_h)$ is well defined whenever $v_h \in X_h$ and $\varphi_h \in Y_h$.

\begin{remark}\label{rah}
If we consider the mapping ${\cal T}_h\: :\,X_h\,\longrightarrow\,Y_h$ defined by
$$ {\cal T}_h v_h=\sum_{\alpha=1}^N\:v_h^\alpha\otimes\sigma_\alpha\,\, \,\,\forall v_h \in X_h,$$
the form $a_h$ can be written on $Y_h$ as
\begin{equation}\label{ah}
a_h\:(v_h,\varphi_h) \,=\,\displaystyle\int_{\Omega}\nabla_H v_h\nabla_H  \varphi_h \,+\, \int_{\Omega}\partial _z ({\cal T}_h(v_h) )\:\partial _z  \varphi_h.
\end{equation}
where we denote $\nabla=(\nabla_H,\partial_z)$ with $\nabla_H=(\partial_1,\cdots,\partial_{d-1})$ the horizontal gradient. 
\end{remark}

\section{Stability and convergence analysis} \label{wellpodness}
In this section we study the well-posedness of problem (\ref{eq:poistransalphadis}) in the sense of Hadamard: We prove that it admits a unique solution that depends continuously on the data.  This will be based upon the Banach-Ne\^cas-Babuska theorem (Cf. ~\cite{ernguermond}). 
\par
We consider the space $X_h$ endowed with the following  $H^1_0(\Omega)$ discrete norm:
\begin{equation}\label{normdisH1}
\begin{array}{l}
\|v_h\|_{X_h}\,=\,\left( h\:\displaystyle\sum_{\alpha=1}^N \| \nabla v_h^\alpha\|^2_{0,\omega} + \displaystyle\:\frac{h}{2}\: \left\|\frac{v_h^{1}}{h/2}\right\|_{0,\omega}^2 
+\sum_{\alpha=1}^{N-1}\:h\: \left\|\frac{v_h^{\alpha+1}-v_h^\alpha}{h}\right\|_{0,\omega}^2 \right .\\ \left.+ \displaystyle\:\frac{h}{2}\: \left\|\frac{v_h^{N}}{h/2}\right\|_{0,\omega}^2 \right )^{1/2}\,=\,\left(\| \nabla_H v_h\|^2_{0,\Omega}\,+\,\|\partial_z ({\cal T}_h v_h)\|^2_{0,\Omega}\right )^{1/2},
\end{array}
\end{equation}
and the space $Y_h$ endowed with the $H_0^1(\Omega)$ norm. It holds the
\begin{lemma}\label{isoXhYh}
The mapping $ {\cal T}_h$ is an isomorphism between normed spaces. Moreover, it holds 
\begin{equation}\label{equivnorm}
\|\nabla ({\cal T}_hv_h)\|_{0,\Omega}\;\leq\;\|v_h\|_{X_h}\;\leq\;\sqrt{3}\:\|\nabla ({\cal T}_h v_h)\|_{0,\Omega},\quad \forall v_h\:\in\:X_h.
\end{equation}

\end{lemma}
\textbf{Proof:} By construction ${\cal T}_h$ is linear and bijective. To prove the equivalence of the norms,  let $\varphi_h \in Y_h$.  Denote $\varphi^\alpha=\varphi_h (z_\alpha)$. Some standard calculations yield
$$ \begin{array}{c}
\|\nabla_H \varphi_h\|^2_{0,\Omega}\,=\,\displaystyle\frac{h}{6}\:\int_\omega \:|\nabla_H \varphi^1|^2\,+\,\sum_{\alpha=1}^{N-1}\:\frac{h}{3}\:\left\{\int_\omega |\nabla_H \varphi^\alpha|^2\,+\,\int_\omega |\nabla_H \varphi^{\alpha+1}|^2\right\}\,+\\[2ex] \displaystyle\,+\,\sum_{\alpha=1}^{N-1}\:\frac{h}{3}\:\int_\omega \nabla_H \varphi^\alpha\cdot \nabla_H \varphi^{\alpha+1}\,+\,\frac{h}{6}\int_\omega \,|\nabla_H \varphi^N|^2.
\end{array}
$$
Also, as 
\begin{equation}\label{derzfunctest}
\partial_z\varphi_h\,=\,\left\{\begin{array}{cll}
\displaystyle \frac{\varphi^{1}}{h/2} &\mbox{in} &\omega\times[z_b,\:z_{1}], \\[2ex] 
\displaystyle\frac{(\varphi^{\alpha+1} - \varphi^\alpha)}{h}\:  &\mbox{in} &\omega\times[z_\alpha,\:z_{\alpha+1}], \quad \forall \, \alpha = 1, 2, . . . , N-1 \\[2ex]
\displaystyle  \frac{-\varphi^{N}}{h/2}&\mbox{in} &\omega\times[z_{N},\:z_{s}]; \\[2ex] 
\end{array}\right.
\end{equation}
it follows
\begin{equation}\label{normphi}
\begin{array}{c}
\|\nabla \varphi_h\|^2_{0,\Omega}\,=\,\displaystyle\frac{h}{2}\:\|\nabla_H \varphi^1\|_{0,\omega}^2\,+\,\sum_{\alpha=2}^{N-1}\:\frac{2h}{3}\:\|\nabla_H \varphi^\alpha\|_{0,\omega}^2\,+\\[2ex] \displaystyle\,+\,\sum_{\alpha=1}^{N-1}\:\frac{h}{3}\:\int_\omega \nabla_H \varphi^\alpha\cdot \nabla_H \varphi^{\alpha+1}\,+\,\frac{h}{2}\|\nabla_H \varphi^N\|_{\omega}^2\,+\\[2ex]
+\, \displaystyle\int_{0,\omega} \frac{(\varphi^1)^2}{h/2}\,+\,\displaystyle\sum_{\alpha=1}^{N-1}\:\int_{\omega} \frac{(\varphi^{\alpha+1}-\varphi^\alpha)^2}{h} \,+\, \displaystyle\int_{\omega}  \frac{(\varphi^{N})^2}{h/2}.
\end{array}
\end{equation}
To obtain the first inequality in \eqref{equivnorm}, (\ref{normphi}) and Young's  inequality yield
$$
\begin{array}{c}
\|\nabla \varphi_h\|^2_{0,\Omega}\,\leq\,\displaystyle\frac{h}{2}\:\|\nabla_H \varphi^1\|_{0,\omega}^2\,+\,\sum_{\alpha=2}^{N-1}\:\frac{2h}{3}\:\|\nabla_H \varphi^\alpha\|_{0,\omega}^2\,+\\[2ex] \displaystyle\,+\,\sum_{\alpha=1}^{N-1}\:\frac{h}{6}\:\left(\|\nabla_H \varphi^\alpha\|_{0,\omega}^2+\|\nabla_H \varphi^{\alpha+1}\|_{0,\omega}^2\right)\,+\,\frac{h}{2}\|\nabla_H \varphi^N\|_{0,\omega}^2\,+\\[2ex]
+\, \displaystyle\int_{\omega} \frac{(\varphi^1)^2}{h/2}\,+\,\displaystyle\sum_{\alpha=1}^{N-1}\:\int_{\omega} \frac{(\varphi^{\alpha+1}-\varphi^\alpha)^2}{h}  \,+\, \displaystyle\int_{\omega}  \frac{(\varphi^{N})^2}{h/2} \,\leq\,\|({\cal T}_h)^{-1}(\varphi_h)\|_{X_h}^2.
\end{array}
$$
Analogously, to obtain the first inequality in \eqref{equivnorm}, again (\ref{normphi}) and  Young's  inequality yield
$$
\begin{array}{c}
\|\nabla \varphi_h\|^2_{0,\Omega}\,\geq\,\displaystyle\frac{h}{2}\:\|\nabla_H \varphi^1\|_{0,\omega}^2\,+\,\sum_{\alpha=2}^{N-1}\:\frac{2h}{3}\:\|\nabla_H \varphi^\alpha\|_{0,\omega}^2\,+\\[2ex] \displaystyle\,-\,\sum_{\alpha=1}^{N-1}\:\frac{h}{6}\:\left(\|\nabla_H \varphi^\alpha\|_{0,\omega}^2+\|\nabla_H \varphi^{\alpha+1}\|_{0,\omega}^2\right)\,+\,\frac{h}{2}\|\nabla_H \varphi^N\|_{0,\omega}^2\,+\\[2ex]
+\, \displaystyle\int_{\omega} \frac{(\varphi^1)^2}{h/2}\,+\,\displaystyle\sum_{\alpha=1}^{N-1}\:\int_{\omega} \frac{(\varphi^{\alpha+1}-\varphi^\alpha)^2}{h}  \,+\, \displaystyle\int_{\omega}  \frac{(\varphi^{N})^2}{h/2} \,\geq\,\frac{1}{3}\:\|({\cal T}_h)^{-1}(\varphi_h)\|_{X_h}^2,
\end{array}
$$
what concludes the proof with $v_h\,=\,({\cal T}_h)^{-1}(\varphi_h)$. \hfill $\square$
\par
\vspace*{.2cm}
Observe that this results justifies the choice of the discrete norm \eqref{normdisH1} if $X_h$ is considered as an external approximation of $H^1_0(\Omega)$. 
\par
The stability of the multilayer problem (\ref{eq:poistransalphadis}) is stated next.
\begin{theorem}\label{wellposed}
The form  $a_h(\cdot,\:\cdot)$ satisfies the following properties:  
\begin{align}
1. & \qquad  \quad a_h\:(v_h,{\cal T}_h(v_h))\,\geq\,\frac{1}{2}\:\|v_h\|_{X_h}\, \|\nabla ({\cal T}_hv_h)\|_{0,\Omega}\,\, \forall\;v_h\in X_h; \label{eq:ahest1} \\
2. & \qquad a_h\:(v_h,\varphi_h)\,\leq\,\sqrt{3}\:\|v_h\|_{X_h}\, \|\nabla \varphi_h\|_{0,\Omega},\,\, \forall\;v_h\in X_h,\;\varphi_h\,\in\, Y_h. \label{eq:ahacot1}
\end{align}
\end{theorem}
\textbf{Proof:}  1. Let $v_h\in X_h$ and $\varphi_h\in Y_h$. It holds
\begin{equation}\label{termH}
\begin{array}{l}
\displaystyle\sum_{\alpha=1}^N (\nabla v_h^\alpha,\nabla  \:\varphi_h)_{\Omega_\alpha}=\displaystyle\sum_{\alpha=1}^{N}\:\int_{z_{\alpha-1/2}}^{z_{\alpha+1/2}}\int_\omega \nabla_H v_h^\alpha \cdot\nabla_H \:\varphi_h
\\[2ex]
=\displaystyle\int_{z_b}^{z_{1}}\int_\omega \, \frac{z-z_{b}}{h/2}\: \nabla_H \varphi^{1}\cdot\nabla_H v_h^{1}\,+\,\displaystyle\int_{z_{N}}^{z_{s}}\int_\omega \:\frac{z_{s}-z}{h/2}\: \nabla_H \varphi^{N} \cdot\nabla_H v_h^{N}
\\[2ex]
+\displaystyle\sum_{\alpha=1}^{N-1}\:\int_{z_\alpha}^{z_{\alpha+1/2}}\int_\omega \,\left( \frac{z_{\alpha+1}-z}{h}\:\nabla_H \varphi^\alpha+\frac{z-z_{\alpha}}{h}\:\nabla_H \varphi^{\alpha+1}\right)\cdot\nabla_H v_h^\alpha
\\[2ex]
+\displaystyle\sum_{\alpha=2}^{N}\:\int_{z_{\alpha-1/2}}^{z_{\alpha}}\int_\omega \,\left( \frac{z_{\alpha}-z}{h}\:\nabla_H \varphi^{\alpha-1}+\frac{z-z_{\alpha-1}}{h}\:\nabla_H \varphi^{\alpha}\right)\cdot\nabla_H v_h^{\alpha}
\\[2ex]
=
\displaystyle\left(\frac{3h}{8}+\frac{h}{4}\right)\:\int_\omega \,\nabla_H v_h^{1}\cdot\nabla_H \varphi^{1}\,+\,\left(\frac{3h}{8}+\frac{h}{4}\right)\:\int_\omega \,\nabla_H v_h^{N}\cdot \nabla_H \varphi^{N}
\\[2ex]
+\displaystyle\displaystyle\frac{h}{8}\:\sum_{\alpha=2}^{N}\:\int_\omega\:\nabla_H v_h^\alpha \cdot \nabla_H  \varphi^{\alpha-1}+\frac{3h}{4}\:\displaystyle\sum_{\alpha=2}^{N-1}\int_\omega \nabla_H v_h^\alpha \cdot \nabla_H \varphi^{\alpha}\\
\displaystyle+\frac{h}{8}\:\displaystyle\sum_{\alpha=1}^{N-1}\:\int_\omega\:\nabla_H v_h^\alpha \cdot\nabla_H \varphi^{\alpha+1}.
\end{array}
\end{equation}
Consequently, we can express $a_h\:(v_h,\varphi_h)$, as 
\begin{equation}\label{ahvhfh}
\begin{array}{l}
a_h\:(v_h,\varphi_h)\,=\,\displaystyle\frac{5h}{8}\:\int_\omega \,\nabla_H\:v_h^{1}\cdot\nabla_H \varphi^{1}+\,\displaystyle\frac{5h}{8}\:\int_\omega \,\nabla_H\:v_h^{N}\cdot\nabla_H \varphi^{N}\,
\\[2ex]
+\displaystyle \,\frac{3h}{4}\:\displaystyle\sum_{\alpha=2}^{N-1}\int_\omega \nabla_H v_h^\alpha\cdot \nabla_H \varphi^{\alpha}\\
\displaystyle +\,\frac{h}{8}\:\sum_{\alpha=1}^{N-1}\:\left(\int_\omega\:\nabla_H v_h^\alpha \cdot\nabla_H \varphi^{\alpha+1}\,+\,\int_\omega\: \nabla_H v_h^{\alpha+1}\cdot\nabla_H \varphi^\alpha\right)\,+\\[2ex]
\displaystyle+ \,\int_{\omega} v_h^{1}\:\frac{\varphi^1}{h/2} +\, \displaystyle\int_{\omega}  v_h^{N}\frac{\varphi^{N}}{h/2}   +\sum_{\alpha=1}^{N-1}\:\int_{\omega} (v_h^{\alpha+1}-v_h^\alpha)\frac{\varphi^{\alpha+1}-\varphi^\alpha}{h} \,   .
\end{array}
\end{equation}
 To obtain \eqref{eq:ahest1}, let $\varphi_h \,=\,{\cal T}_h (v_h)$. Using  (\ref{ahvhfh}), we have
 \begin{equation}\label{ahvhTvh}
\begin{array}{l}
a_h\:(v_h,\varphi_h)\,=\,\displaystyle \displaystyle\frac{5h}{8}\:\|\nabla_H v_h^1\|_{0,\omega}^2\,+\frac{5h}{8}\:\|\nabla_H v_h^N\|_{0,\omega}^2
\\[2ex]
\displaystyle+\,\frac{3h}{4}\:\displaystyle\sum_{\alpha=2}^{N-1}\|\nabla_H v_h^\alpha\|_{0,\omega}^2+\displaystyle\frac{h}{4}\:\sum_{\alpha=1}^{N-1}\:\int_\omega\:\nabla_H v_h^\alpha \cdot\nabla_H v_h^{\alpha+1}\,\\
\displaystyle+\int_{\omega} \frac{(v_h^1)^2}{h/2}\,+
 \int_{\omega}  \frac{(v_h^{N})^2}{h/2} +\displaystyle\sum_{\alpha=1}^{N-1}\:\int_{\omega} \frac{(v_h^{\alpha+1}-v_h^\alpha)^2}{h}   .
\end{array}
\end{equation}
By Young's inequality,
$$
\int_\omega\:\nabla_H v_h^\alpha \cdot\nabla_H v_h^{\alpha+1}\ge -\frac{1}{2}\,\left (\|\nabla_H v_h^\alpha\|_{0,\omega}^2\,+\,\| \nabla_H v_h^{\alpha+1}\|_{0,\omega}^2 \right ).
$$
Then, from (\ref{ahvhTvh}), it follows
\begin{equation}\label{estiminfsup}
\begin{array}{c}
a_h\:(v_h,\varphi_h)\,\geq \displaystyle
\frac{h}{2}\:\sum_{\alpha=1}^{N}\:\|\nabla_H v_h^\alpha\|_{0,\omega}^2 +\displaystyle\int_{\omega} \frac{(v_h^1)^2}{h/2}+ \displaystyle\int_{\omega}  \frac{(v_h^{N})^2}{h/2}
\\[2ex]
+\,\displaystyle\sum_{\alpha=1}^{N-1}\int_{\omega} \frac{(v_h^{\alpha+1}-v_h^\alpha)^2}{h}\geq\,\displaystyle\frac{1}{2}\:\|v_h\|^2_{X_h}
\,\geq\,\frac{1}{2}\:\|v_h\|_{X_h}\,\|\nabla \varphi_h\|_{0,\Omega}.
\end{array}
\end{equation}
 \par\noindent
2. Let  $v_h\in X_h$ and $\varphi_h\in Y_h$, starting from (\ref{ahvhfh}) and using (\ref{equivnorm}) and  Holder's inequality, we have

 $$
\begin{array}{l}
a_h\:(v_h,\varphi_h)\,\leq\,\displaystyle\frac{5h}{8}\: \|\nabla_H\:v_h^{1}\|_{0,\omega}\:\|\nabla_H \varphi^{1}\|_{0,\omega}\,+\,\displaystyle\frac{5h}{8}\:\|\nabla_H\:v_h^{N}\|_{0,\omega}\:\|\nabla_H \varphi^{N}\|_{0,\omega}\,+\\[2ex]
+\,\displaystyle\sum_{\alpha=2}^{N-1}\: \frac{3h}{4}\: \|\nabla_H\:v_h^{\alpha}\|_{0,\omega}\: \:\|\nabla_H \varphi^{\alpha}\|_{0,\omega}
\,+\,\displaystyle\sum_{\alpha=1}^{N-1}\:\frac{h}{8}\: \|\nabla_H\:v_h^{\alpha}\|_{0,\omega}\: \:\|\nabla_H \varphi^{\alpha+1}\|_{0,\omega}\,\\[2ex]
+\displaystyle\sum_{\alpha=1}^{N-1}\:\,\frac{h}{8}\: \|\nabla_H\:v_h^{\alpha+1}\|_{0,\omega}\:\|\nabla_H \varphi^{\alpha}\|_{0,\omega}\,+\,  \left(\int_{\Gamma_b}\frac{(v_h^{1})^2}{h/2}\right)^{1/2}\:\left(\int_{\Gamma_b} \frac{(\varphi^1)^2}{h/2}\right)^{1/2}\,+\\[2ex]
+\displaystyle\sum_{\alpha=1}^{N-1}\left(\int_{\Gamma_{\alpha+\frac{1}{2}}} \frac{(v_h^{\alpha+1}-v_h^\alpha)^2}{h}\right)^{1/2}\left(\int_{\Gamma_{\alpha+\frac{1}{2}}} \frac{(\varphi^{\alpha+1}-\varphi^\alpha)^2}{h}\right)^{1/2} \,+\,\\[2ex]\,+\displaystyle\left(\int_{\Gamma_s}\frac{(v_h^{N})^2}{h/2}\right)^{1/2}\left(\int_{\Gamma_s} \frac{(\varphi^N)^2}{h/2}\right)^{1/2}
\, \leq \end{array}
$$
$$
\begin{array}{l}
\leq\,\left\{\displaystyle\frac{5h}{8}\: \|\nabla_H\:v_h^{1}\|_{0,\omega}^2\,+\,\displaystyle\frac{5h}{8}\:\|\nabla_H\:v_h^{N}\|_{0,\omega}^2\,+\,\displaystyle\sum_{\alpha=2}^{N-1}\: \frac{3h}{4}\: \|\nabla_H\:v_h^{\alpha}\|_{0,\omega}^2
\,+\right.\\[2ex]
+\,\displaystyle\sum_{\alpha=1}^{N-1}\:\frac{h}{8}\: \|\nabla_H\:v_h^{\alpha}\|_{0,\omega}^2\,
+\,\displaystyle\sum_{\alpha=1}^{N-1}\:\,\frac{h}{8}\: \|\nabla_H\:v_h^{\alpha+1}\|_{0,\omega}^2\,+\sum_{\alpha=1}^{N-1}\int_{\Gamma_{\alpha+\frac{1}{2}}} \frac{(v_h^{\alpha+1}-v_h^\alpha)^2}{h}\,+\\[2ex]
\left.+\,\displaystyle \int_{\Gamma_b}\frac{(v_h^{1})^2}{h/2}
\, +\,\int_{\Gamma_s}\frac{(v_h^{N})^2}{h/2}\right\}^{1/2}\cdot\left\{\displaystyle\frac{5h}{8}\: \|\nabla_H\:\varphi^{1}\|_{0,\omega}^2\,+\,\displaystyle\frac{5h}{8}\:\|\nabla_H\:\varphi^{N}\|_{0,\omega}^2\,+\right.\\[2ex]
+\,\displaystyle\sum_{\alpha=2}^{N-1}\: \frac{3h}{4}\: \|\nabla_H\:\varphi^{\alpha}\|_{0,\omega}^2
\,+\,\displaystyle\sum_{\alpha=1}^{N-1}\:\frac{h}{8}\: \|\nabla_H\:\varphi^{\alpha}\|_{0,\omega}^2\,
+\,\displaystyle\sum_{\alpha=1}^{N-1}\:\,\frac{h}{8}\: \|\nabla_H\:\varphi^{\alpha+1}\|_{0,\omega}^2\,+\\[2ex]
\left.\displaystyle+\,\sum_{\alpha=1}^{N-1}\int_{\Gamma_{\alpha+\frac{1}{2}}} \frac{(\varphi^{\alpha+1}-\varphi^\alpha)^2}{h}\,+\,\displaystyle \int_{\Gamma_b}\frac{(\varphi^{1})^2}{h/2}\, +\,\int_{\Gamma_s}\frac{(\varphi^{N})^2}{h/2}\right\}^{1/2}\,\leq \\[2ex]
\leq\,\|v_h\|_{X_h}\:\|({\cal T}_h)^{-1}(\varphi_h)\|_{X_h}\leq\,\sqrt{3}\:\|v_h\|_{X_h}\:\|\nabla \varphi_h\|_{0,\Omega}.
\end{array}
$$
\hfill $\square$
\par
\vspace*{.2cm}
As a conclusion we deduce 
\begin{corollary}\label{stable0}
The form $a_h$ satisfies the inf-sup condition
\begin{equation}\label{eq:17}
\displaystyle \displaystyle \inf_{v_h\in X_h}\sup_{\varphi_h\in Y_h}\frac{a_h\:(v_h,\varphi_h)}{\|\nabla \varphi_h\|_{0,\Omega}\:\|v_h\|_{X_h}}\,\geq\;\frac{1}{2}.
\end{equation}
\end{corollary}
\textbf{Proof:} It is a direct consequence of estimate \eqref{eq:ahest1}.\hfill $\square$
\par \vspace*{.4cm}
This implies the well posedness of problem (\ref{eq:poistransalphadis}):
\begin{corollary}\label{stable}
The multilayer problem (\ref{eq:poistransalphadis}) admits a unique solution $v_h \in X_h$ that satisfies the estimate 
\begin{equation} \label{eq:estvh}
\|v_h\|_{X_h} \,\leq\,2\:\|f\|_{H^{-1}(\Omega)}.
\end{equation}
\end{corollary}
\textbf{Proof:} The form $a_h$ is stable by Theorem \ref{wellposed}, and satisfies the inf-sup condition \eqref{eq:17}. Also, \eqref{eq:ahest1} implies that if $a_h\:(v_h,\varphi_h) = 0$ for all $v_h\,\in\,X_h, $ for some $\varphi_h\:\in\:Y_h$, then necessarily $\varphi_h\;=\,0.$
Consequently the hypotheses of the Banach-Necas-Babuska Theorem hold (Cf. ~\cite{ernguermond}). This ensures that problem (\ref{eq:poistransalphadis}) admits a unique solution that depends continuously on the data $f$. The constant $2$ in estimate \eqref{eq:estvh} also follows from \eqref{eq:ahest1}. \hfill $\square$
%
\subsection{Convergence analysis}\label{ConvergeSection}
In this section  we prove a convergence result for general solutions of problem \eqref{fv:poispb}, so as optimal order error estimates for smooth solutions. For that, let us construct an interpolation operator on $X_h$ for functions defined on $\Omega$. Let ${\cal P}_h$ be the prismatic grid of $\Omega$ constructed by vertical displacements of the grid ${\cal T}_k$, located at the nodes $z_b,\:z_{1}$,\: $z_{2} \:\cdots,\: z_{N-1},\: z_{N},\:z_s$. The geometric elements of ${\cal P}_h$ are constructed as $P=K \times I_j$, where $I_1=[z_{b},z_{1}]$, $I_j=[z_{j-1},z_{j}]$, $j=2,\cdots,N$; $I_{N+1}=[z_{N},z_{s}]$ and $K \in {\cal T}_k$ . Note that $Y_h$ is the prismatic finite element space defined by
$$
 Y_h=\{ w \in C^0(\Omega) \,|\,w_{|_P} \in R_l(K) \times P_1(I_j), \,\forall \, P=K \times I_j \in {\cal P}_h \,\}\,\cap\,H_0^1\:(\Omega).
$$
If $\{x_i\}_{i\in {\cal I}}$ is the set of Lagrange interpolation nodes for $V_{0,k}(\omega)$, for some set of indices ${\cal I}$, then the set of Lagrange interpolation nodes for $Y_h$ is \\ $\{r_{i,j}=(x_i, z_{j}) \,| \,i\in{\cal I},\, j=1,\cdots, N\,\,  \}$. 
Following Bernardi et al. \cite{bernardi}, we introduce a nodal Lagrange interpolation operator $\Pi_k: L^2(\Omega) \mapsto Y_h$, of the form
$$
(\Pi_k w )(x,z) =\sum_{i\in {\cal I}} \sum_{j= 1}^{N}\bar{w}_{i,j} \, \alpha_{i,j}(x,z),
$$
where the $\alpha_{i,j}$ are the nodal Lagrange interpolation basis functions on $V_{0,k}(\omega)$ corresponding interpolation nodes $r_{i,j}$, and the $\bar{w}_{i,j}$ are  averaged value of $w$ on an element to which the node $r_{i,j}$ belongs.

The operator $\Pi_k$ satisfies the standard finite element interpolation error estimates,
\begin{equation} \label{eq:errestpik1}
\|w - \Pi_k w\|_{0,\Omega} \le \hat C\, (h+k^l)\, \|w\|_{l,\Omega},\,\mbox{ if } w \in H^l (\Omega),
\end{equation}
\begin{equation} \label{eq:errestpik2}
\|w - \Pi_k w\|_{1,\Omega} \le \hat C\, (h+k^l)\, \|w\|_{l+1,\Omega},\,\mbox{ if } w \in H^{l+1}(\Omega).
\end{equation}
We shall assume that the operator $\Pi_k$ also is defined componentwise for vector functions, without change of notations.  We next define an interpolation operator
$\Pi_{h}: \, H^1_0(\Omega) \,\longrightarrow\:X_h\quad \mbox{that transforms } v\,\in\, H^1_0(\Omega)$ into $\Pi_{h} v \,\in\, X_h,$
defined by $(\Pi_{h} v)|_{\Omega_\alpha}\,=\,v^\alpha_k\:\in\:V_{0,k}(\omega),\quad\forall\; \alpha=1,2,\cdots,N$, with
\begin{equation}\label{defPih}
v^\alpha_k\:(x)\,=\,\displaystyle\frac{1}{h}\:\int_{z_{\alpha-\frac{1}{2}}}^{z_{\alpha+\frac{1}{2}}}\:\Pi_k v(x,z)\:dz \quad x\;\in \;w\;\mbox{a. e.}.
\end{equation}
We at first obtain the optimal order error estimate
 \begin{theorem}\label{converge1}
Assume that the weak solution of (\ref{fv:poispb}) satis\-fies $v \in H^{l+1}(\Omega)\cap H^1_0(\Omega)$, with $l\ge 1$. Then, the solution $v_h\,\in\, X_h$ of the discretization (\ref{eq:poistransalphadis}) verifies 
\begin{equation}\label{esterr}
\|\Pi_{h} v\,-\,v_h\|_{X_h}\,\leq\:C\:\|v\|_{l+1,\Omega}\,(h+k^l)
 \end{equation}
for some constant $C>0$ a constant.
\end{theorem}
\textbf{Proof:} Define the consistency error $\varepsilon_h\,\in\,Y_h'$ by
$$ <\varepsilon_h,\varphi_h>\,=\,a_h\:(\Pi_{h} v\,-\,v_h,\varphi_h),\quad \forall \varphi_h \in Y_h.$$
It is sufficient show that, there exits a constant $C'>0$ such that
\begin{equation}\label{estimerr}
|<\varepsilon_h,\varphi_h>|\,\leq\:C'\:\|v\|_{l+1,\Omega} \:\|\nabla \varphi_h\|_{0,\Omega}\,(h+k^l),\quad \forall\;\varphi_h\,\in\, Y_h.
\end{equation}
Indeed, by the inf-sup condition, we have 
$$
\begin{array}{c}
\displaystyle \frac{1}{2}\:\|\Pi_{h} v\,-\,v_h\|_{X_h}\,\leq\,\sup_{\varphi_h\,\in\, Y_h}\:\frac{a_h\:(\Pi_{h} v\,-\,v_h,\varphi_h)}{\|\nabla \varphi_h\|_{0,\Omega}}\,=\\[2ex]\,=\,\displaystyle\sup_{\varphi_h\,\in\, Y_h}\:\frac{<\varepsilon_h,\varphi_h>}{\|\nabla \varphi_h\|_{0,\Omega}}
\leq\:C'\:\|v\|_{l+1,\Omega} \,(h+k^l),
\end{array}
$$
and we obtain the error estimate with $C=2\,C'$.

To prove~(\ref{estimerr}), notice that $a_h\:(v_h,\varphi_h)=L(\varphi_h)=a\:(v,\varphi_h)$, and then
$$
\begin{array}{c}
<\varepsilon_h,\varphi_h>\,=\,a_h\:(\Pi_{h} v\,-\,v_h,\varphi_h)\,=\,a_h\:(\Pi_{h} v,\varphi_h)\,-\,a\:(v,\varphi_h)\,= \\[2ex]
=\,\displaystyle\sum_{\alpha=1}^N \:\int_{\Omega_\alpha}\nabla (v^\alpha-v)\nabla  \varphi_h \,+ \,\displaystyle\int_{\Gamma_b} v^1 \partial _z\varphi_h\,+\, \displaystyle\sum_{\alpha=1}^{N-1}\:\int_{\Gamma_{\alpha+\frac{1}{2}}} (v^{\alpha+1}-v^\alpha) \partial _z\varphi_h \, +\\[2ex]    
- \displaystyle\int_{\Gamma_s} v^N \partial _z\varphi_h\, =\,\displaystyle\sum_{\alpha=1}^N \:\int_{\Omega_\alpha}\nabla_H (v^\alpha-v)\nabla_H  \varphi_h \,-\, \sum_{\alpha=1}^N \:\int_{\Omega_\alpha}\partial _z v \:\partial _z  \varphi_h \,+\\[2ex]
\displaystyle+\,\int_{\omega}\int_{z_{b}}^{z_{1}} \:\frac{v^{1}}{h/2} \partial _z\varphi_h + \displaystyle\sum_{\alpha=1}^{N-1}\:\int_{\omega}\int_{z_{\alpha}}^{z_{\alpha+1}}\:\frac{v^{\alpha+1}-v^\alpha}{h} \partial _z\varphi_h + \int_{\omega}\int_{z_{N}}^{z_{s}}\: \frac{-v^{N}}{h/2} \partial _z\varphi_h
=\end{array}
$$
$$
\begin{array}{c}
=\displaystyle\sum_{\alpha=1}^N \:\int_{\Omega_\alpha}\nabla_H (\Pi_{h} v\,-\,v)\nabla_H  \varphi_h \,-\, \sum_{\alpha=1}^N \:\int_{\Omega_\alpha}\partial _z v \:\partial _z  \varphi_h \,+\,\int_{\omega}\int_{z_{b}}^{z_{1}} \partial _z ({\cal T}_h(\Pi_{h} v)) \:\partial _z  \varphi_h+
\\[2ex]
\displaystyle+\, \sum_{\alpha=1}^{N-1}\:\int_{\omega}\int_{z_{\alpha}}^{z_{\alpha+1}}\: \partial _z ({\cal T}_h(\Pi_{h} v))\:\partial _z  \varphi_h\:   \,+\,\displaystyle\int_{\omega}\int_{z_{N}}^{z_{s}}\:\partial _z ({\cal T}_h(\Pi_{h} v)) \:\partial _z  \varphi_h .
\end{array}
$$
Then 
\begin{equation}\label{errcons}
<\varepsilon_h,\varphi_h>\,=\,\displaystyle\int_{\Omega}\nabla_H (\Pi_{h} v\,-\,v)\cdot\nabla_H  \varphi_h \,+\, \int_{\Omega}\partial _z ({\cal T}_h(\Pi_{h} v) \,-\,v)\:\partial _z  \varphi_h.
\end{equation}
To estimate the first term in (\ref{errcons}), 
we have
$$
\displaystyle\left|\int_{\Omega}\nabla_H (\Pi_{h} v\,-\,v)\cdot\nabla_H  \varphi_h\right|\,\leq\, \|\nabla_H (\Pi_{h} v\,-\,v)\|_{0,\Omega} \|\nabla_H  \varphi_h\|_{0,\Omega}.
$$
Now we add $\pm (\nabla_H v)(x,\hat z)$ to the integral, use Holder's and Minkowski's inequalities and apply that if $v \in H^2(\Omega)$, it holds
$$\displaystyle
(\nabla_H v)(x,\hat z)\,=\,(\nabla_H v)(x,z)\,+\,\int_{\hat z}^{z}\:\partial _z (\nabla_H v)(x,s)ds$$  $\mbox{for  }z,\:\hat z\in(z_{\alpha-\frac{1}{2}}, z_{\alpha+\frac{1}{2}})$, a. e. for $x \in \omega$.
Then, 
$$
\begin{array}{l}
\displaystyle\|\nabla_H (\Pi_{h} v)\,-\,\nabla_H v\|_{0,\Omega}\,= \\[2ex] =\,\displaystyle
\left(\sum_{\alpha=1}^N \:\int_{\Omega_\alpha} \left|\frac{1}{h}\:\int_{z_{\alpha-\frac{1}{2}}}^{z_{\alpha+\frac{1}{2}}}\:(\nabla_H (\Pi_{k} v)(x,\hat z)\,-\,(\nabla_H  v)(x,z))\:d\hat z\right|^2 dx dz\right)^{1/2}\leq
\\[2ex] \leq\,\displaystyle \left(\frac{1}{h}\:\sum_{\alpha=1}^N \:\int_{\Omega_\alpha} \int_{z_{\alpha-\frac{1}{2}}}^{z_{\alpha+\frac{1}{2}}}\:\left|(\nabla_H (\Pi_{k} v)(x,\hat z)\,-\,(\nabla_H  v)(x,z))\right|^2 \:d\hat z dx dz\right)^{1/2}\,\le
\end{array}
$$

\begin{equation}\label{estimH}
\begin{array}{l}
\leq\displaystyle\left(\frac{1}{h}\:\sum_{\alpha=1}^N \:\int_{\Omega_\alpha} \int_{z_{\alpha-\frac{1}{2}}}^{z_{\alpha+\frac{1}{2}}}\:\left|(\nabla_H (\Pi_{k} v)(x,\hat z)\,-\,(\nabla_H  v)(x,\hat z))\right|^2 \:d\hat z dx dz\right)^{1/2}+
\\[2ex]
+\,\displaystyle\left(\frac{1}{h}\:\sum_{\alpha=1}^N\int_{\Omega_\alpha} \:\int_{z_{\alpha-\frac{1}{2}}}^{z_{\alpha+\frac{1}{2}}}\:\left|((\nabla_H v)(x,\hat z)\,-\,(\nabla_H  v)(x,z))\right|^2 \:d\hat z dx dz\right)^{1/2}\,=
\\[2ex]
=\,\displaystyle\left(\sum_{\alpha=1}^N \:\int_{\Omega_\alpha}\:\left|(\nabla_H (\Pi_{k} v)(x,\hat z)\,-\,(\nabla_H  v)(x,\hat z))\right|^2 \: dx d\hat z\right)^{1/2}\,+
\\[2ex]
+\,\displaystyle\left(\frac{1}{h}\:\sum_{\alpha=1}^N\:\int_{\Omega_\alpha} \:\int_{z_{\alpha-\frac{1}{2}}}^{z_{\alpha+\frac{1}{2}}}\:\left|\int_{\hat z}^{z}\:\partial _z (\nabla_H v)(x,s)ds\right|^2 \:d\hat z dx dz\right)^{1/2}\,\leq
\\[2ex]
\displaystyle\leq\,\|\nabla_H (\Pi_{k} v\,-\, v)\|_{0,\Omega}\,+\,\displaystyle h\:\|\partial _z (\nabla_H v)\|_{0,\Omega} \,\le
\\[2ex]
\displaystyle\leq\,\hat C\, (h+k^l)\, \|v\|_{l+1,\Omega}\,+\,\displaystyle h\:\|\partial _z (\nabla_H v)\|_{0,\Omega} .
\end{array}
\end{equation}
where the last inequality follows from (\ref{eq:errestpik2}). 

To estimate the term in $\partial _z$ in (\ref{errcons}), we have
$$
\left|\int_{\Omega}\partial _z ({\cal T}_h(\Pi_{h} v) \,-\,v)\:\partial _z  \varphi_h\right|\,\leq\,\displaystyle\|\partial _z ({\cal T}_h(\Pi_{h} v) \,-v)\|_{0,\Omega}\,\|\partial _z  \varphi_h\|_{0,\Omega}. 
$$
Observe that the composed operator ${\cal T}_h \circ \Pi_h$ is a Cl\'ement interpolation operator (cf. \cite{clement}), so as $v \in H^{2}(\Omega)$, it holds
\begin{equation}\label{estimdz}
\|\partial _z ({\cal T}_h(\Pi_{h} v) \,-\,\partial _z v\|_{0,\Omega}\,\leq\,\hat C'\:\|v\|_{2,\Omega}\:h.
\end{equation}
 Then, from (\ref{estimH}) and (\ref{estimdz}), we obtain (\ref{estimerr}) with $C'=\hat C' + \hat C +1$.\hfill $\Box$
\par
\vspace*{.4cm}
\begin{remark}
The overall order of the discretization \eqref{eq:poistransalphadis} is then one. Thus, to minimize the number of degrees of freedom the best choice is $l=1$ and $k=h$.
\end{remark}
We now can deduce the general convergence result,
\begin{corollary}\label{coroconv}
Let $v \in H^1_0(\Omega)$ the solution of problem (\ref{fv:poispb}). Then it holds
\begin{equation}\label{gralconv}
\lim_{h \to 0} \|\Pi_h v - v_h\|_{X_h}=0
\end{equation}
where $v_h$ is the solution of the multilayer discretization (\ref{eq:poistransalphadis}). 
\end{corollary}
\textbf{Proof:}  As in the proof of Theorem~\ref{converge1}, it holds
\begin{eqnarray}
\frac{1}{2}\|\Pi_h v - v_h\|_{X_h} &\le &\sup_{\varphi_h \in Y_h}\displaystyle\frac{<\varepsilon_h,\varphi_h>}{\|\nabla \varphi_h\|_{0,\Omega}} \nonumber\\
& \le& \label{eserdos}
\|\nabla_H (\Pi_{h} v \,-\, v)\|_{0,\Omega}
+\|\partial _z (({\cal T}_h\circ\Pi_{h}) v \,-\, v)\|_{0,\Omega}.
\end{eqnarray}
Operators $\Pi_h$ and ${\cal T}_h\circ\Pi_{h}$ are stable, in the sense that
$$
\|\nabla_H (\Pi_{h} w)\|_{0,\Omega} \le C\, \|w\|_{1,\Omega},\,\, \forall w \in H^1(\Omega),
$$
$$
\|\partial _z (({\cal T}_h\circ\Pi_{h}) w)\|_{0,\Omega} \le C\, \|\partial _z w\|_{0,\Omega},\,\, \forall w \in H^1(\Omega).
$$
Indeed, is straightforward to prove that $\|\nabla_H (\Pi_{h} w)\|_{0,\Omega} \le \|\Pi_k w\|_{1,\Omega}$ and then the first estimate follows as it because $\Pi_k$ is stable in $H^1(\Omega)$ norm. The second one is a standard property of the Cl\'ement interpolation operator (cf. \cite{clement}). Then a standard argument, using the density of $ {\cal D}(\Omega)$ in $H^1_0(\Omega)$, proves that the r. h. s. of \eqref{eserdos} tends to zero as $h \to 0$. \hfill $\Box$

 \section{Extension to domains with non-flat upper boundary} \label{se:nonflat}
In this section we extend the multilayer discretization of the Poisson problem (\ref{eq:poispb0}) to domains with non-flat upper boundary, that arise in  geophysical flows.  Concretely, we consider domains $\hat \Omega$ that can be obtained as vertical transformations of the domain  $\Omega\,=\,\omega\:\times \:(0,1)$ (that we call in this section \lq\lq reference" domain) by a change of variables of the form:
\begin{equation}\label{Kchange}
\hat z \,= \,\displaystyle z\,\eta(x)
, \quad \hat x \,= \,x, 
\end{equation}
with $x\:\in\:\omega\subset\,\mathbb{R}^2$, $z\:\in\:(0,1)$ and where $\eta$ is an smooth and strictly positive function defined in $\omega$ (the surface equations are thus $z=\eta(x,y)$). We assume that $\eta\:\in\:W^{1,\infty}(\omega)$ is a function such that
\begin{equation}\label{condeta}
\displaystyle\eta(x)\,\geq\,\eta_0\,>\,0\; \mbox{and} \; \|\nabla_H \eta\|_{\infty,\Omega}\,<\,1, \; \forall\;x\in\omega\;\mbox{uniformly},
\end{equation}
for some constant $\eta_0$. Note that he restriction in \eqref{condeta} on the gradient of $\eta$ is consistent with the need of having $\eta >0$ on $\omega$. 

We intend to solve the Poisson problem in $\hat \Omega$:
\begin{equation}\label{eq:poispb2}
\left\{\begin{array} {ccccc}
-\Delta \hat v &= & \hat f & \mbox{in} & \hat \Omega, \\ [2ex]
\hat v & = & 0 & \mbox{on} & \partial \hat \Omega.
\end{array}\right.
\end{equation}

To do it, we use the change of variables (\ref{Kchange}), to transform problem  (\ref{eq:poispb2}) in the  following elliptic problem in the reference domain $\Omega$ in variational form:  \\ Find $v \in H_0^1(\Omega)$  such that
\begin{equation}\label{fv:KLap}
a_K\:(v,\varphi) \,=  \, L\:(\varphi),\qquad \forall \;\varphi \,\in\,H_0^1(\Omega), 
\end{equation}
with $a_K\:(v,\varphi) \,= \,\displaystyle\int_{\Omega}\:\eta\:(\nabla  \varphi)^T\: K\: \nabla v $ and $L\:(\varphi)\,=\,<f , \varphi>_{\Omega}$, where $f\:(x,z)\:=\:\hat f (x, \hat z)$ and $K$ is a symmetric and positive definite matrix given by 
$$K \,=\, \left(\begin{array} {cc} K_{HH}  & K_{H3}  \\ [2ex]
K_{3H}  & K_{33} \end{array}\right)$$
where $K_{HH}\,=\,\displaystyle Id_H$ with $Id_H$ the $2\times 2$ identity matrix, $K_{H3}\,=\,\displaystyle-\:\frac{z}{\eta}\:\nabla_H \eta$, $K_{3H}\,=\,K_{H3}^T$ and $K_{33}\,=\, \displaystyle\frac{1 }{\eta^2}\,+\,|K_{H3}|^2$

Here our objective is to deduce a multilayer approach of problem (\ref{eq:poispb2}). 

\subsection{Multilayer discretization}

To obtain our new multilayer system, we consider the same vertical decomposition of domain $\Omega$ and the same definitions and notations as those used in Section 2. 

Similarly to the definition of the form $a_h$ by \eqref{ah}, we propose the following multilayer discretization of problem (\ref{eq:poispb2}): 
\begin{equation}\label{fv:KLaptransalphadis}
\left\{\begin{array}{l}
\mbox{Find}\quad v_h\,\in\, X_h, \quad \mbox{such that}
\\[2ex]
a_{K,h}\:(v_h,\varphi_h) \,=\, L\:(\varphi_h), \qquad \forall \;\varphi_h\, \in \,Y_h,
\end{array}\right. 
\end{equation}
where, the bilinear form $a_{K,h}\:(\cdot,\cdot)$ is given by
\begin{equation}\label{aKh}
\begin{array} {l}
a_{K,h}\:(v_h,\varphi_h) \,=\,\displaystyle\int_{\Omega}\:\eta\left[\:(\nabla_H v_h)^T\: \nabla_H  \varphi_h\,  + \, \partial _z ({\cal T}_h(v_h))   K_{3H} \nabla_H \varphi_h  \right]\,+
\\[2ex]
+\,\displaystyle\int_{\Omega}\:\eta\left[(\nabla_H v_h)^T  K_{3H}^T \partial _z  \varphi_h \, +\,  \partial _z ({\cal T}_h(v_h))   K_{33} \partial _z \varphi_h \right].
\end{array}
\end{equation}

\subsection{Analysis of the multilayer discretization}
We now study some properties of the multilayer problem (\ref{fv:KLaptransalphadis}).
The stability of  the multilayer discretization is achieved, under some restrictions on the gradient of function $\eta$, as follows:
\begin{theorem}\label{Kwellposed}
We assume that function $\eta$ satisfies (\ref{condeta}). Then, the bilinear form  $a_{K,h}(\cdot,\:\cdot)$ satisfies the following properties:
\begin{enumerate}
\item  There exists a constant $C_{1,K}>0$ that depends on $\eta$, such that 
\begin{equation}\label{equivKinfsup}
a_{K,h}\:(v_h,{\cal T}_h(v_h))\,\geq\,C_{1,K}\:\|\nabla {\cal T}_h(v_h)\|_{0,\Omega}\:\|v_h\|_{X_h},\qquad \forall\;v_h\in X_h.
\end{equation}

\item There exists a constant $C_{2,K}>0$ that
depends on $\eta$, such that
\begin{equation}\label{eq:aKhacot1}
a_{K,h}\:(v_h,\varphi_h)\,\leq\,C_{2,K}\:\|\nabla \varphi_h\|_{0,\Omega}\:\|v_h\|_{X_h},\; \forall\;v_h\in X_h,\; \forall\;\varphi_h\,\in\, Y_h.
\end{equation}
\end{enumerate}
\end{theorem}

\textbf{Proof:}
\par\noindent
1. To obtain (\ref{equivKinfsup}), let $v_h\in X_h$ and $\varphi_h \,=\,{\cal T}_h (v_h)$, then from (\ref{aKh}) we have
$$
\begin{array} {l}
a_{K,h}\:(v_h,{\cal T}_h (v_h)) \,= \,\displaystyle\int_{\Omega}\:\eta\left[\:(\nabla_H v_h)^T\: \nabla_H  ({\cal T}_h (v_h))\, +\, K_{33} |\partial _z {\cal T}_h (v_h)|^2 \,  \right]\,+
\\[2ex]
+\,\displaystyle\int_{\Omega}\:\eta\left[\partial _z {\cal T}_h(v_h)   K_{3H} \nabla_H\: ( {\cal T}_h (v_h)+ v_h) \right].
\end{array}
$$

Now, using that $\eta\,>\,0$,  Young's  inequality $\left(\displaystyle a\:b  \geq -\frac{1}{2}\:(\epsilon a^2+\frac{1}{\epsilon}b^2)\right)$ for the last integral and the estimates within the proof of  (\ref{equivnorm}) in Lemma~\ref{isoXhYh} for the integrals with $\nabla_H v_h$ or $\nabla_H  ({\cal T}_h (v_h))$, we obtain

$$
\begin{array} {l}
a_{K,h}\:(v_h,{\cal T}_h (v_h))\,\geq\,\displaystyle\int_{\Omega}\:\eta\left[\:\frac{1}{2}|\nabla_H  v_h|^2\, +\, K_{33} |\partial _z {\cal T}_h (v_h)|^2 \,  \right] \,+
\\[2ex]
-\,\displaystyle\frac{\epsilon}{2}\:\int_{\Omega}\:\eta   |K_{3H}|^2 |\partial _z {\cal T}_h(v_h)|^2\,-\,\displaystyle\frac{1}{\epsilon}\:\int_{\Omega}\:\eta  
 \left(|\nabla_H\: ( {\cal T}_h (v_h))|^2\,+\,|\nabla_H\: v_h|^2\right) \,\geq
\\[2ex]
\geq\,\displaystyle\left(\frac{1}{2}-\frac{2}{\epsilon}\right)\:\int_{\Omega}\:\eta|\nabla_H  v_h|^2\,+\,\int_{\Omega}\:\eta   \left(K_{33}-\frac{\epsilon}{2}|K_{3H}|^2\right) |\partial _z {\cal T}_h(v_H)|^2\,.
 \end{array}
$$
Then, if we take $\varepsilon\,=\,4\:(1+\beta)$ with $\beta\,>\,0$ and we use that $ \|\nabla_H \eta\|_{\infty,\Omega}\,<\,1$, we have
$$
\begin{array} {l}
a_{K,h}\:(v_h,{\cal T}_h (v_h))\,\geq\,\displaystyle \,\int_{\Omega}\:\left(\frac{1-(1+2\beta)z^2(|\nabla_H\eta|^2)}{\eta}\right) |\partial _z {\cal T}_h(v_h)|^2\,+
\\[2ex]
\displaystyle+\,\frac{\beta}{2(1+\beta)}\:\int_{\Omega}\:\eta|\nabla_H  v_h|^2\,\geq\,C_{1,K}\:\|v_h\|_{X_h}^2\,\geq\,C_{1,K}\:\|\nabla({\cal T}_h (v_h))\|_{0,\Omega}\:\|v_h\|_{X_h} .
 \end{array}
$$
with $C_{1,K}\,=\,\displaystyle\min\left\{\frac{\beta\:\eta_0}{2(1+\beta)},\frac{1-(1+2\beta)\|\nabla_H \eta\|_{\infty,\Omega}^2}{\|\eta\|_{\infty,\Omega}}\right\}$, where for the last inequality we use  (\ref{equivnorm}).
\par\noindent
2. Let $v_h\,\in\, X_h$ and $\varphi_h\,\in\, Y_h$, then as $\eta\,>\,0$, (\ref{equivnorm}) and HolderÕs inequality, we have
$$
\begin{array}{c}
a_{K,h}\:(v_h,\varphi_h)\,\leq\,\displaystyle\left(\int_{\Omega}\:\eta|\nabla_H  v_h|^2\right)^{1/2}\:\left(\int_{\Omega}\:\eta|\nabla_H  \varphi_h|^2\right)^{1/2}\,+
\\[2ex]
+\,\displaystyle\left(\int_{\Omega}\:\eta|K_{3H}|^2 |\partial _z {\cal T}_h(v_h)|^2\right)^{1/2}\:\left(\int_{\Omega}\:\eta|\nabla_H  \varphi_h|^2\right)^{1/2}\,+\,\displaystyle\left(\int_{\Omega}\:\eta|\nabla_H  v_h|^2\right)^{1/2}\:\cdot
\\[2ex]
\displaystyle\cdot\:\left(\int_{\Omega}\:\eta|K_{3H}|^2 |\partial _z \varphi_h|^2\right)^{1/2}+
\left(\int_{\Omega}\:\eta K_{33} |\partial _z {\cal T}_h(v_H)|^2\right)^{1/2}\:\left(\int_{\Omega}\:\eta K_{33} |\partial _z \varphi_h|^2\right)^{1/2}\leq
\\[2ex]
\displaystyle\leq \left\{2\int_{\Omega}\:\eta|\nabla_H  v_h|^2+\int_{\Omega}\:\eta (|K_{3H}|^2+K_{33}) |\partial _z {\cal T}_h(v_h)|^2\right\}^{1/2}\cdot
\\[2ex]
\displaystyle\cdot \left\{2\int_{\Omega}\:\eta|\nabla_H  \varphi_h|^2+\int_{\Omega}\:\eta (|K_{3H}|^2+K_{33}) |\partial _z \varphi_h|^2\right\}^{1/2}\leq
\\[2ex]
\displaystyle\leq C_{K} \|v_h\|_{X_h} \|({\cal T}_h)^{-1}(\varphi_h)\|_{X_h}\,\leq\,\sqrt{3}\:C_K\:\|v_h\|_{X_h}\:\|\nabla \varphi_h\|_{0,\Omega}.
 \end{array}
$$
Thus, we obtain (\ref{eq:aKhacot1}) with $C_{2,K}=\sqrt{3}\:C_{K}$ where $C_{K}=\displaystyle\min\left\{2\|\eta\|_{\infty,\Omega},\frac{1+2\|\nabla_H \eta\|_{\infty,\Omega}^2}{\eta_0}\right\}$.

As consequence of this theorem we deduce the following well posedness result for the multilayer problem (\ref{fv:KLaptransalphadis}):
\begin{corollary}\label{Kstable}
We assume that function $\eta$ satisfies (\ref{condeta}), then the multilayer problem (\ref{fv:KLaptransalphadis}) admits a unique solution $v_h \in X_h$ that satisfies the estimate 
$$ 
\|v_h\|_{X_h} \,\leq\,C_{3,K}\:\|f\|_{H^{-1}(\Omega)}
$$
where  $C_{3,K}\,>\,0$ is a constant that depends on $\eta$.
\end{corollary}

The convergence of the multilayer discretization (\ref{fv:KLaptransalphadis}) is stated as follows: 
\begin{theorem}\label{Kestima}
We assume that $\eta$ verifies (\ref{condeta}) and that the weak solution of (\ref{fv:KLap}) satis\-fies $v \in H^{l+1}(\Omega)\cap H_0^{1}(\Omega)$, with $l\ge 1$.  Then, the solution $v_h\,\in\, X_h$ of the discrete problem (\ref{fv:KLaptransalphadis}) verifies the estimate
$$
\|\Pi_{h} v\,-\,v_h\|_{X_h}\,\leq\,(h+k^l)\:C_{4,K}\:\|v\|_{l+1,\Omega}
 $$
with $C_{4,K}$ a constant that depends on $\eta$ and $\Pi_{h}$ the operator defined by (\ref{defPih}).

Moreover, if $v \in H^1_0(\Omega)$, then, it holds
\begin{equation}\label{gralKconv}
\lim_{h \to 0} \|\Pi_h v - v_h\|_{X_h}=0
\end{equation}
\end{theorem}

The proof of this result is an extension of that of Theorem~\ref{converge1} and Corollary~\ref{coroconv} by standard techniques, we omit it for brevity.

\section{Neumann boundary conditions} \label{se:neumann}
\begin{figure}
  \begin{center}
\makebox[\linewidth]{\includegraphics[height=.24\linewidth]{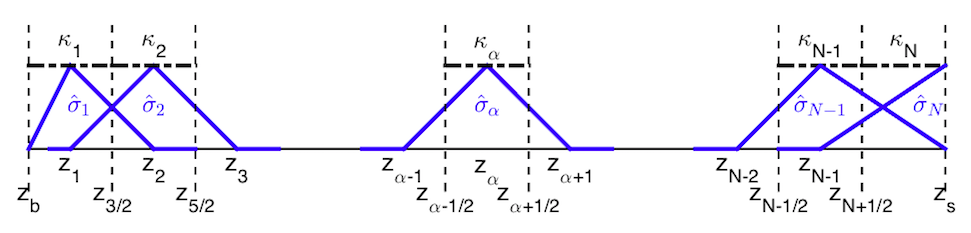}}
  \caption{Basis function on $X_h$ and $\hat Y_h$ \label{fig:basisNeu}}
  \end{center}
\end{figure}

In this section we extend the multilayer discretization to the Poisson problem with Dirichlet boundary conditions  (\ref{eq:poispb0}) to the following Poisson problem with Dirichlet-Neumann boundary conditions:
\begin{equation}\label{eq:neumann}
\left\{\begin{array} {ccccc}
-\Delta v &= & f & \mbox{in} & \Omega, \\ [2ex]
v & = & 0 & \mbox{on} & \Gamma_b\cup \Gamma_l, \\ [2ex]
\partial_n v & = & g & \mbox{on} &  \Gamma_s .
\end{array}\right.
\end{equation}
The variational form of this problem is: Given $f \in H^{-1}(\Omega)$ and $g \in H^{-1/2}\:(\Gamma_s)$, find $v \in H_{bl}^1(\Omega)$  such that
\begin{equation}\label{fv:neumann}
a\:(v,\hat \varphi) \,=  \,\hat L\:(\hat\varphi),\qquad \forall \;\hat\varphi \,\in\,H_{bl}^1(\Omega), 
\end{equation}
with 
\begin{equation} \label{formaNeu}
a\:(v,\hat \varphi) \,=\,\displaystyle\int_\Omega \nabla v \:\nabla \hat \varphi\, d\Omega, \quad \mbox{and}\quad \hat L\:(\hat \varphi)\,=\,\displaystyle<f , \hat\varphi>_\Omega+\int_{\Gamma_s} g \hat \varphi,
\end{equation}
where $H_{bl}^1(\Omega)\displaystyle=\{v \in H^1(\Omega)\, /\; v_{|_{\Gamma_b\cup \Gamma_l}}=0\} $.
Here we consider the same space of semi-discrete solutions in the vertical direction as in Section~\ref{se:multilayer}, that we  endow with the following  $H_{bl}^1(\Omega)$ discrete norm:
\begin{equation}\label{normdisHbl}
\begin{array}{rcl}
|||v_h|||_{X_h}\,&=&\,\left ( h\:\displaystyle\sum_{\alpha=1}^N \| \nabla v_h^\alpha\|^2_{0,\omega} + \displaystyle\sum_{\alpha=1}^{N-2}\:h\: \left\|\frac{v_h^{\alpha+1}-v_h^\alpha}{h}\right\|_{0,\omega}^2 \right .\\
&&\quad \left .+\displaystyle\:\frac{h}{2}\: \left\|\frac{v_h^{1}}{h/2}\right\|_{0,\omega}^2 
+ \displaystyle\:\frac{3h}{2}\: \left\|\frac{v_h^{N}-v_h^{N-1}}{3h/2}\right\|_{0,\omega}^2 \right )^{1/2}.
\end{array}
\end{equation}
This new discrete norm is motivated by the following new semi-discrete test functions space, 
$$\hat Y_h\,=\,\displaystyle\left\{\hat\varphi_h\,\in\,L^2(\Omega)\, /\;\hat \varphi_h\,=\,\sum_{\alpha=1}^N\:\hat\varphi^\alpha\otimes\hat\sigma_\alpha,\;\mbox{with}\;\hat\varphi^\alpha\:\in\:V_{0,k}(\omega)\right\},$$
where $\hat \sigma_i = \sigma_i$, for all $i = 1, . . . , N-2$, and  $\hat \sigma_{N-1}$, $\hat \sigma_{N}$ are piecewise affine 1D functions on the intervals $ [z_{N-2},z_{N-1}]$, $[z_{N-1},z_s]$, respectively, that satisfy $\hat \sigma_{N-1}(z_j)\,=\delta_{(N-1)j}$ and $\hat \sigma_N (z_j)\,=0$, for all $j = 1, 2, . . . , N-1 $; $\hat \sigma_{N-1}(z_b) =\hat  \sigma_{N-1}(z_s) =\hat  \sigma_N(z_b) = 0$ and $\hat \sigma_N(z_s)=1$ as we see in Figure~\ref{fig:basisNeu}. Note that $\hat Y_h$ is a subspace of $H_{bl}^1(\Omega)$ such that $dim (\hat Y_h)=dim(X_h)$ and the interval $[z_{N-1},z_s]$ has length $3h/2$.

We look for an approximation of the solution of (\ref{fv:neumann}) of the form \\ $v \simeq v_h \,+\, g_h$ with $g_h=g_k\:(z-z_s)\:\kappa_N$  and $v_h$ is the solution of the following multilayer discretization of (\ref{fv:neumann}):
\begin{equation}\label{eq:neumanndis}
\left\{\begin{array}{l}
\mbox{Find}\quad v_h\,\in\, X_h, \quad \mbox{such that}
\\[2ex]
\hat a_h\:(v_h,\hat\varphi_h) \,=\,\hat L_h\:(\hat\varphi_h), \qquad \forall \;\hat\varphi_h\, \in \,\hat Y_h.
\end{array}\right. 
\end{equation}
where $\hat a_h$ is a the bilinear form defined by
$$
\begin{array}{rcl}
\hat a_h\:(v_h,\hat \varphi_h) &=&\displaystyle\int_{\Omega}\nabla_H v_{h} \nabla_H  \hat \varphi_h \,+\, \int_{\Omega}\partial _z (\hat{\cal T}_h(v_{h}))\:\partial _z  \hat\varphi_h;
\end{array}
$$
with $\hat  {\cal T}_h$ the mapping defined by $\hat {\cal T}_h   v_h=\displaystyle\sum_{\alpha=1}^N\:  v_h^\alpha\otimes\hat \sigma_\alpha \:\in\:\hat Y_h$, for $  v_h \in  X_h$, and 
$\hat L_h\:(\hat\varphi_h)\,=\,\displaystyle<f , \hat\varphi_h>_\Omega+\int_{\Gamma_s} g \hat\varphi_h - \hat a_h\:(g_h,\hat\varphi_h) $. We assume in particular $g \in H^1_0(\omega)$ in order to have this last term well defined.

\begin{remark}\label{rahNeu} The basis functions of spaces $Y_h$ and $\hat Y_h$ only differ on the layers $\Omega_{N-1}$ and $\Omega_N$. Then the form $\hat a_h$ only differs from the form $a_h$ defined in \eqref{eq:formah} in the integrals on $\Omega_{N-1}\cup \Omega_N$. Then a slight modification of the proof of Lemma \ref{isoXhYh} yields 
\begin{equation}\label{tauhNeu}
\frac{4}{5}\|\nabla (\hat{\cal T}_hv_h)\|_{0,\Omega}\;\leq\;|||v_h|||_{X_h}\;\leq\;4\:\|\nabla (\hat {\cal T}_hv_h)\|_{0,\Omega},\,\forall v_h\:\in\;X_h.
\end{equation}
Consequently, the mapping $\hat  {\cal T}_h$ is too an isomorphism between the normed spaces $X_h$ and $\hat Y_h$. This justifies the choice of the new norm on space $X_h$.
\end{remark}

\subsection{Well-posedness and convergence analysis}\label{PropertiesNeumannSection}
In this section we first study, analogously as we see in Section \ref{wellpodness}, the well-posedness of the multilayer problem \eqref{eq:neumanndis}:
\begin{theorem}\label{NeumannStable}
We assume that $g\in H_0^1(\omega)$, then the multilayer problem (\ref{eq:neumanndis}) admits a unique solution $v_h \in X_h$ that satisfies the estimate 
\begin{equation} \label{eq:estvhNeu}
|||v_h|||_{X_h} \,\leq\,\frac{5}{2}\:\left(\|f\|_{H^{-1}(\Omega)}\,+\, (\sqrt{L}+2\sqrt{h})\:\|g\|_{0,\omega}\,+\, h\sqrt{h}\:\|\nabla_H g\|_{0,\omega}\right).
\end{equation}
\end{theorem}
\textbf{Proof:} (Sketch) Reasoning as in the proof of Theorem \ref{wellposed} for the multilayer problem (\ref{eq:poistransalphadis}), we deduce that  the form $\hat a_h$ is stable and satisfies the inf-sup condition
\begin{equation}\label{InfsupNeumann}
\displaystyle \displaystyle \inf_{v_h\in X_h}\sup_{\hat \varphi_h\in\hat Y_h}\frac{\hat a_h\:(v_h,\hat \varphi_h)}{\|\nabla\hat \varphi_h\|_{0,\Omega}\:|||v_h|||_{X_h}}\,\geq\;\frac{2}{5}.
\end{equation}
This new inf-sup condition is proved similarly to \eqref{eq:17}, as a consequence of estimates \eqref{tauhNeu} and taking into account the differences between the forms $a_h$ and $\hat a_h$ mentioned in Observation~\ref{rahNeu}.

Then problem (\ref{eq:neumanndis}) admits a unique solution that depends continuously on the data $f$ and $g$. The estimate  \eqref{eq:estvhNeu} follows from \eqref{InfsupNeumann}, and standard estimates of $\hat L_h(\hat \varphi_h)$. 
\hfill $\square$

Finally  we prove, analogously again to the case of the Poisson problem \eqref{fv:poispb}, the convergence of the discrete solution of the multilayer system (\ref{eq:neumanndis})  to the weak solution of the variational problem \eqref{fv:neumann}, so as optimal order error estimates for smooth solutions. For that, let us define the operator  $\hat\Pi_{h}$ defined by $\hat \Pi_{h} v \,=\,\Pi_{h} v + g_h$, for $v \in H^1_{bl}(\Omega)$, with  $\Pi_{h}$ defined by \eqref{defPih}. Then we have:
 \begin{theorem}\label{convergeneumann}
Assume that the weak solution of (\ref{fv:neumann}) satis\-fies $v \in H^{l+1}(\Omega)\cap H^1_{bl}(\Omega)$, with $l \ge 1$. Then, the solution $v_h\,\in\, X_h$ of the discretization (\ref{eq:neumanndis}) verifies 
\begin{equation}\label{esterrneumann}
|||\hat\Pi_{h} v\,-\,v_h|||_{X_h}\,\leq\:C\:((h+k^l)\|v\|_{l+1,\Omega}\,+\,\sqrt{h}\:\|g\|_{0,\omega}+ h\:\sqrt{h}\:\|\nabla_H g\|_{0,\omega}) \end{equation}
for some constant $C>0$ a constant.

Moreover, if the solution $v$ of problem (\ref{fv:neumann}) is only in $H^1_{bl}(\Omega)$, then it holds
\begin{equation}\label{gralconvneumann}
\lim_{h \to 0} |||\hat\Pi_h v -v_h|||_{X_h}=0.
\end{equation} 
\end{theorem}

\section{Numerical results} \label{se:tests}
In this section we present some numerical 3D experiments, to analyze the computing time reduction obtained by the parallel computation of each multilayer discretization, so as to test the error estimates. We have used the FreeFem++ software (cf. \cite{freefem++}).

Concretely, for the Poisson problem~(\ref{eq:poispb0}), we solve the multilayer problem ~(\ref{eq:poistransalphadis}) by an iterative procedure. Taking the test functions $\varphi_h=\varphi^\alpha\otimes\sigma_\alpha$,  $1\leq\alpha\leq N$ in  ~(\ref{eq:poistransalphadis}), this problem is equivalent to the following linear system, with unknowns $v_h^{\alpha}, \, \alpha=1,\cdots, N \in V_{0,k}(\omega)$:

$$
\begin{array}{l}
\displaystyle\int_\omega \:\left(\frac{5h}{8}\:\nabla_H\:v_h^{1} \,+\, \frac{h}{8}\:\nabla_H\:v_h^{2}\right)\cdot\nabla_H \:\varphi^{1} \,+\,\int_{\omega}\:\left(\frac{3}{h}\: v_h^{1}\,-\, \frac{1}{h}\:v_h^{2}\right)\:\varphi^1 \,=
\\[3ex]
\displaystyle=\,L(\varphi^1\otimes\sigma_1), \quad \forall \varphi^1 \in V_{0,k}(\omega);
\end{array}
$$
$$
\begin{array}{l}
\displaystyle\int_\omega \:\left(\frac{h}{8}\nabla_H v_h^{\alpha-1}\,+\,\frac{3h}{4}\:\nabla_H v_h^\alpha\,+\, \frac{h}{8}\nabla_H v_h^{\alpha+1}\right)\cdot \nabla_H \varphi^{\alpha}\,+
\\[2ex]
+\displaystyle \int_\omega \left(-\,\frac{1}{h}\:v_h^{\alpha-1}\,+\,\frac{2}{h}\:v_h^\alpha\,-\,\frac{1}{h}\:v_h^{\alpha+1}\right)\: \varphi^{\alpha}\,=\,L(\varphi^\alpha\otimes\sigma_\alpha), \quad \forall \varphi^\alpha \in V_{0,k}(\omega),\\ \alpha=2,\cdots, N-1;
\\[2ex]
\displaystyle\int_\omega \:\left(\frac{h}{8}\:\nabla_H\:v_h^{N-1} \,+\, \frac{5h}{8}\:\nabla_H\:v_h^{N}\right)\cdot\nabla_H \:\varphi^{N} \,+\,\int_{\omega}\:\left(-\frac{1}{h}\: v_h^{N-1}\,+\, \frac{3}{h}\:v_h^{N}\right)\:\varphi^N \,=
\\[3ex]
\displaystyle=\,L(\varphi^N\otimes\sigma_N),  \quad \forall \varphi^N \in V_{0,k}(\omega).
\end{array}
$$

We solve this linear system through a block-Jacobi iterative algorithm by layers, that leads to solve at each iteration the following  sequence of 2D horizontal problems : 
\begin{equation}\label{Jacobi}
\begin{array}{l}
\displaystyle\int_\omega \:\frac{5h}{8}\:\nabla_H\:v_h^{1,k+1}\cdot\nabla_H \varphi^{1} \,+\, \frac{3}{h}\:\int_{\omega} v_h^{1,k+1}\varphi^1\,= -\frac{h}{8}\:\int_\omega \:\nabla_H\:v_h^{2,k}\cdot\nabla_H \varphi^{1}\,+
\\[2ex]
\displaystyle  +\, \frac{1}{h}\:\int_{\omega}v_h^{2,k}\:\varphi^1 \,+\,\frac{5h}{8}\:\int_\omega\:f^{1}\:\varphi^1+\frac{h}{8}\:\int_\omega\:f^{2}\:\varphi^1,  \quad \forall \varphi^1 \in V_{0,k}(\omega);
\\
\\[2ex]
\displaystyle\frac{3h}{4}\:\int_\omega \nabla_H v_h^{\alpha,k+1} \cdot \nabla_H \varphi^{\alpha}\,+\,\frac{2}{h}\:\int_\omega \:v_h^{\alpha,k+1}\: \varphi^{\alpha}\,=\,
\\[2ex]
\displaystyle - \:\frac{h}{8}\int_\omega\:\left(\nabla_H v_h^{\alpha-1,k}\,+\,\nabla_H v_h^{\alpha+1,k}\right)\cdot \nabla_H \varphi^{\alpha}\,
+\displaystyle \frac{1}{h}\:\int_\omega \left(v_h^{\alpha-1,k}\,+\,v_h^{\alpha+1,k}\right)\: \varphi^{\alpha}
\\[2ex]
+\displaystyle\int_\omega\:\left(\frac{h}{8}\:f^{\alpha-1}+\frac{3h}{4}\:f^\alpha+\frac{h}{8}\:f^{\alpha+1}\right) \varphi^{\alpha}, 
 \quad \forall \varphi^\alpha \in V_{0,k}(\omega)\; 2\leq\alpha\leq N-1;
\\
\\[2ex]
\displaystyle\int_\omega \:\frac{5h}{8}\:\nabla_H\:v_h^{N,k+1}\cdot\nabla_H \varphi^{N} \,+\, \frac{3}{h}\:\int_{\omega} v_h^{N,k+1}\varphi^N\,= -\frac{h}{8}\:\int_\omega \:\nabla_H\:v_h^{N-1,k}\cdot\nabla_H \varphi^{N}\,+
\\[2ex]
\displaystyle  +\, \frac{1}{h}\:\int_{\omega}v_h^{N-1,k}\:\varphi^1 \,+\,\frac{5h}{8}\:\int_\omega\:f^{N}\:\varphi^N+\frac{h}{8}\:\int_\omega\:f^{N-1}\:\varphi^N,  \quad \forall \varphi^N \in V_{0,k}(\omega),
\end{array}
\end{equation} 
with $\displaystyle f^\alpha\,=\,\int_\Omega\:f\:\sigma_\alpha$,  $\alpha=1,\cdots, N$.

This block-Jacobi algorithm is solved using an affine parallel GMRES algorithm (built-in in FreeFem++), where each layer-wise 2D problem is solved in a different processor using basically the following function:
\begin{lstlisting}
real[int] func Jzero(real[int] & U)	// jacobi iteration
    {ISendRecv(U,uh); 
      bh=rhs;
      w = Mi*uh[0][] ; bh += w ;
      w = Mi*uh[2][] ; bh += w ;
      w = A^-1*bh;
      U-=w;   
      return U;}
\end{lstlisting}
where, for example, to solve the 2D problem on layer $\alpha$ with $\alpha=2,\cdots, N-1$, \verb+rhs+ is the contribution of the integrals in $f$ in \eqref{Jacobi}, \verb+Mi*uh[0][]+ are the integrals in $v_h^{\alpha-1}$ in \eqref{Jacobi} and \verb+Mi*uh[2][]+ are the integrals in $v_h^{\alpha+1}$.

The numerical tests have been run on  Altix UV 2000 with 32  CPUs Intel Xeon 64 bits EvyBridge E4650 with  10 core with a total of 2048 Gb RAM , using the same number of processors as of layers.

\subsection*{Test 1: Homogeneous Dirichlet boundary conditions}
In this numerical test, we consider a smooth exact solution 
$$v(x,y,z)\,=\,4\:z^2\:(1-z)\:\sin(\pi x^2)\:\sin(\pi y)$$ 
of ~(\ref{eq:poispb0}) in $\Omega\,=\,]0,1[^3$ with homogeneous Dirichlet boundary conditions.  We compare the solution obtained by the block-Jacobi algorithm (\ref{Jacobi}), with a piecewise affine finite element space $V_{0,k}(\omega)$ constructed on a structured horizontal grid of $\omega\,=\,]0,1[^2$, versus the direct solution of the global 3D problem (\ref{eq:poispb0}) by the standard Galerkin method with piecewise affine finite elements on a structured tetrahedral grid, with the same degrees of freedom as $X_h$ (that is, number of layers $\times \, dim(X_h)$). 

Concretely, we at first compare the CPU time to obtain the sequential solution by a Conjugate Gradient Solver (built-in in FreeFem++ too) versus the mean CPU time required by the processors to solve the block-Jacobi algorithm (\ref{Jacobi}). 

For that, we compute two tests. Assume that the horizontal grid size is $k=1/NH$ for some integer number $NH \ge 1$. In the first test we set $NZ=NH$ and progressively increase these  numbers from $NZ=NH= 10$ to $NZ=NH=100$, with step $10$. In the second one, we fix the horizontal grid size to the value $NH=1/40$ and progressively increase the  number of layers from $NZ= 10$ to $NZ=100$, with step $10$. 

\begin{figure}
  \begin{center}
  \makebox[\linewidth]{\includegraphics[height=.5\linewidth]{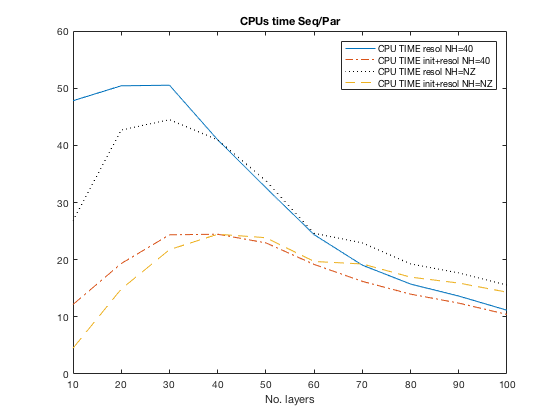}}
  \caption{Test 1: Ratio sequential versus  parallel CPU times.
   \label{fig:solhom}}
  \end{center}
\end{figure}

Figure \ref{fig:solhom} shows for each of these two tests two ratios, on the one hand, the ratio between the CPU time to solve the sequential 3D problem by the Conjugate Gradient Solver (CG)  versus the mean CPU time used by processors to build the matrix and the second member of the 2D system that we solve in each layer/processor and solve it by the affine parallel GMRES algorithm (lines with legend \verb|"CPU TIME init+resol"|). On the other hand, the ratio between the CPU time to solve the sequential problem by the CG Solver  versus the mean CPU time used by processors only to solve the 2D system in each layer/processor by the affine parallel GMRES algorithm (lines with legend \verb|"CPU TIME resol"|)

We notice that we obtain similar results for both tests. As the number of layers increases, there is a large gain of computing time, up to an optimal rate around 30/40 layers, that deteriorates as $NZ$ increases from this value.

Furthermore, we test the error estimates. For that, Table 1 shows the obtained results for the relative errors in the $L^2(\Omega)$ and the $H^1(\Omega)$ discrete norms  introduced in (\ref{normdisH1}), for $NZ= 10, 20, 40, 80$ and $k=1/NZ$, thus with the same grid size in the horizontal and vertical directions. Also, Table 1 displays the estimated orders of convergence, computed from these errors and the number of iterations used by the affine parallel GMRES algorithm to solve the 2D system in each layer/processor.

In all test we denote by $\|v_h\|_{1,h}\,=\,\|v_h\|_{X_h}$ and $\|v_h\|_{0,h}\,=\,\left (\displaystyle\:h\:\sum_{\alpha=1}^{N}\:\|v_h^\alpha\|_{0,\omega}^2\,\right )^{1/2}$, and we compute $ord_{k,h}=\displaystyle \log_2\left(\frac{\|e_h\|_{k,h}}{\|e_{h/2}\|_{k,h/2}}\right)$ for $k=0$, 1. 
\begin{table}
\begin{center}
\begin{tabular}{|c|c|c|c|c|}
\hline
N & 10 & 20 & 40 & 80 
\\ \hline
$\|e_h\|_{0,h}$ & 0.0311414 & 0.0077288 & 0.00192875 & 0.000481821
\\ \hline
$\|e_h\|_{1,h}$ & 0.162031 & 0.080199 & 0.0399963 & 0.0199851
\\ \hline
GMRES iters & 21 &   39 & 89 & 252
\\ \hline
$ord_{0,h}$ & 2.0105 & 2.0026 & 2.0011 & -
\\ \hline
$ord_{1,h}$ & 1.0146 & 1.0037 & 1.009 & -
\\ \hline
\end{tabular}
\caption{Test 1, relative errors in discrete $L^2(\Omega)$ and $H^1(\Omega)$ norms and estimated convergence orders.}
\end{center}
\end{table}

We obtain first order convergence in the $H^1(\Omega)$ discrete norm, in full agreement with our theoretical expectations. We also obtain second order convergence in the $L^2(\Omega)$ discrete norm, as arises for Galerkin finite element solutions of regular elliptic problems. 

We obtain quite similar convergence orders if we use unstructured horizontal meshes, that we not display for brevity.

\section*{Test 2: Domains with non-flat surface}
Now we present a numerical test, where we consider a computational domain $\hat \Omega$  that is obtained from the reference domain $\Omega\,=\,]0,1[^3$ by the change of variables (\ref{Kchange}) with the function $\eta$ given by 
$$\eta(x,y)\,=\,1\,+\,\epsilon\:\sin (2\pi(x+y)).$$ 
This function satisfies (\ref{condeta}) if $\epsilon\,<\,\frac{1}{2\sqrt{2}\pi}$. 

\begin{figure}[htp]
  \begin{center}
  \makebox[\linewidth]{\includegraphics[height=.5\linewidth]{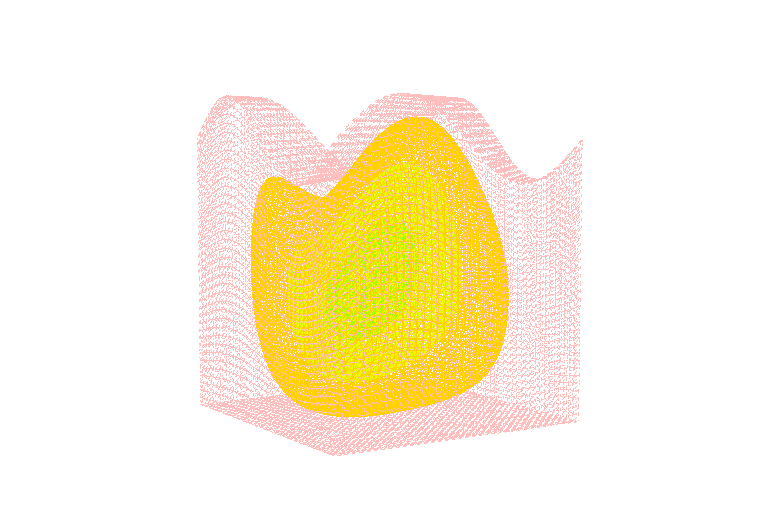}}
  \caption{Test 2. Solution of problem (\ref{eq:poispb2})  with 36 layers.
   \label{fig:Ksol}}
  \end{center}
\end{figure}

In this test we compare an exact solution of the Laplace problem~(\ref{eq:poispb2}) in  $\hat \Omega$ with the solution obtained by our multilayer discretization (\ref{fv:KLaptransalphadis}) in the unit cube $\Omega$. Concretely, 
we consider the function 
$$u(x,y,z)\,=\,64\:\sin\left(\frac{\pi z}{ \eta(x,y)}\right)\:x^2\:(1-x)\:y\:(1- y)^2$$ 
that is the exact solution for the Laplace problem~(\ref{eq:poispb2}) in  $\hat \Omega$, and  we consider the previous function $\eta$ with $\epsilon =0.15$. Figure~\ref{fig:Ksol} shows the global 3D Galerkin solution of (\ref{eq:poispb2}) with $h=1/36$.

\begin{figure}
  \begin{center}
 \makebox[\linewidth]{\includegraphics[height=.5\linewidth]{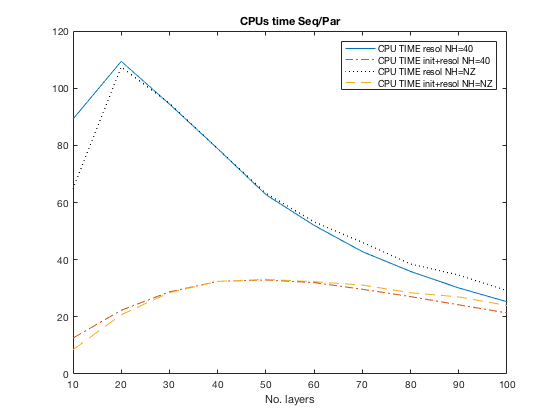}}
   \caption{Test 2: Ratio sequential versus  parallel CPU times.
   \label{fig:Kratio}}
  \end{center}
\end{figure}

Figure~\ref{fig:Kratio} shows the same two CPU ratios for each of the two tests shown in Figure  \ref{fig:solhom}

We observe a similar behavior as in Test 1, there is an increasing computational gain as the number of processors increases up to 40/50 layers (lines with legend \verb|"CPU TIME init+resol"|), that further progressively deteriorates.The gain is somewhat higher than in the previous test, possibly due to the larger computational cost to obtain the 3D sequential global  solution in the computational domain $\hat \Omega$.

Also,  Table~2 display the norms of the relative errors, estimated orders of convergence, as well as the number of iterations used by GMRES, for $NZ= 10, 20, 40, 80$. We observe a quite similar behavior as in Test 1, we recover the theoretical first order in the $H^1(\Omega)$ discrete norm, and second order in the $L^2(\Omega)$ discrete norm.
\begin{table}[ht]
\begin{center}
\begin{tabular}{|c|c|c|c|c|}
\hline
N & 10 & 20 & 40 & 80 
\\ \hline
$\|e_h\|_{0,h}$ & 0.0284935 & 0.00677639 & 0.00158417 & 0.000390097
\\ \hline
$\|e_h\|_{1,h}$ & 0.162359 & 0.0786338 & 0.0381604 & 0.0189835
\\ \hline
GMRES iters & 24 &  46 & 117 & 355
\\ \hline
$ord_{0,h}$ & 2.0720 & 2.0968 & 2.0218 & -
\\ \hline
$ord_{1,h}$ & 1.0460 & 1.0431 & 1.0073 & -
\\ \hline
\end{tabular}
\caption{Test 2, relative errors in discrete $L^2(\Omega)$ and $H^1(\Omega)$ norms and estimated convergence orders.}
\end{center}
\end{table}
Here we do the tests for structured horizontal meshes, but we again obtain quite similar results if we use unstructured horizontal meshes.

\section*{Test 3: Dirichlet-Neumann boundary conditions}
\begin{figure}
  \begin{center}
  \makebox[\linewidth]{\includegraphics[height=.5\linewidth]{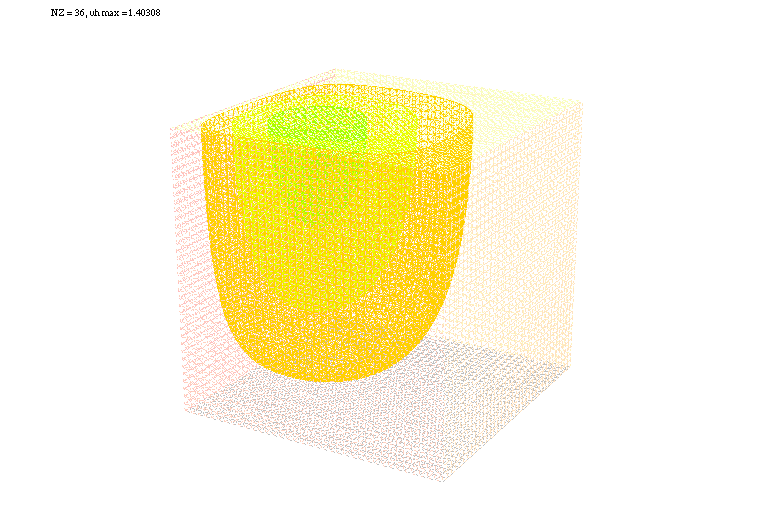}}
  \caption{Test 3. Solution of problem~(\ref{eq:neumann})  with 36 layers.
   \label{fig:solNeumann}}
  \end{center}
\end{figure}

In this third test, we consider an exact solution of ~(\ref{eq:neumann}) in $\Omega$ with Neumann boundary conditions on $\Gamma_s$ and homogeneous Dirichlet boundary conditions on $\Gamma_b\cup\Gamma_l$: 
$$v(x,y,z)\,=\,64\:\sin\left(\frac{\pi z}{ 2}\right)\:x^2\:(1-x)\:y\:(1- y)^2.$$ 
Figure~\ref{fig:solNeumann} shows the global 3D Galerkin solution of ~(\ref{eq:neumann}) with $h=1/36$.

\begin{figure}
  \begin{center}
 \makebox[\linewidth]{\includegraphics[height=.5\linewidth]{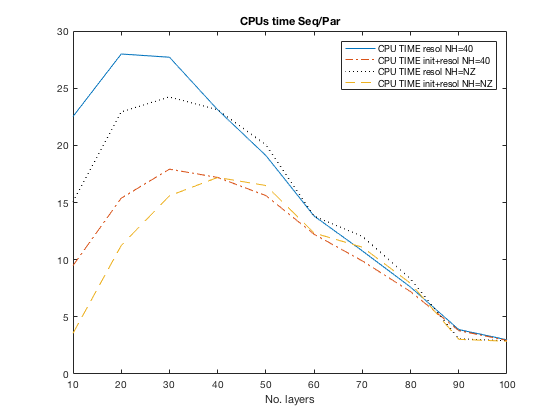}}
   \caption{Test 3. Ratio multilayer versus Galerkin 3D CPU times.
   \label{fig:ratioNeumann}}
  \end{center}
\end{figure}

Figure~\ref{fig:ratioNeumann} shows shows the same two CPU ratios for each of the two tests as in Figure  \ref{fig:solhom} and Figure~\ref{fig:Kratio}

We observe a similar behavior, there is an increasing computational gain as the number of processors increases to around 30/40 layers, that further progressively deteriorates. The gain is somewhat smaller as for the Test 1 due to additional calculations required in the upper layers.

Also,  Table~3 displays the norms of the relative errors, estimated orders of convergence and number of GMRES iterations, for structured horizontal meshes, again for $NZ= 10, 20, 40, 80$. We observe a quite similar behavior to the Test 1 that qualitatively are the same for unstructured horizontal meshes.
\begin{table}[ht]
\begin{center}
\begin{tabular}{|c|c|c|c|c|}
\hline
N & 10 & 20 & 40 & 80 
\\ \hline
$\|e_h\|_{0,h}$ & 0.0423782 & 0.0104868 & 0.00261555 & 0.000654319
\\ \hline
$\|e_h\|_{1,h}$ & 0.21214 & 0.10533 & 0.0526301 & 0.0263241
\\ \hline
GMRES iters & 22 &  44 & 122 & 635
\\ \hline
$ord_{0,h}$ & 2.0147 & 2.0034 & 1.9990 & -
\\ \hline
$ord_{1,h}$ & 1.0101 & 1.0010 & 0.9995 & -
\\ \hline
\end{tabular}
\caption{Test 3, relative errors in discrete $L^2(\Omega)$ and $H^1(\Omega)$ norms and estimated convergence orders.}
\end{center}
\end{table}

\section{Conclusions}
We have studied in this paper a multilayer discretization of second order elliptic problems, aimed at providing reliable multilayer discretizations of shallow fluid flow problems with diffusive effects. This is a Petrov-Galerkin discretization in which the trial functions are piecewise constant per horizontal layers, while the trial functions are continuous piecewise linear, on a vertically shifted grid.

We have introduced the discretization for the Poisson problem with Dirichlet boundary conditions in cylindric domains, and extended it to Neumann boundary conditions in cylindric domains, so as to domains with variable surface, that appear in free-surface fluid flow problems. We have proved the well posedness and optimal error order estimates for these three discretizations in natural norms, based upon specific inf-sup conditions.

We have performed numerical tests with parallel computing of the solution by a block-Jacobi algorithm, based upon the multilayer structure of the discretization, for academic problems with smooth solutions. We recover the theoretical optimal error order convergence, and observe a high increase of the speed of computations for a moderate number of processors. This confirms the interest of applying the technique introduced to multilayer discretizations of fluid flow problems.

In further steps we shall improve the parallelization procedure to obtain a good scaling of the CPU computing time for large number of processors. Also, we will apply the discretization introduced to convection-diffusion problems, and then to multilayer discretizations of Navier-Stokes and related equations. We will also extend the discretization to higher-order approximations of the unknowns. These works are now in progress.


\section*{Acknowledgements}
This research was partially supported by \lq\lq Proyecto de Excelencia de la Junta de Andaluc\'{\i}a"$\;$  under grant P12-FQM-454.

\section*{References}

\bibliography{mybibfile}

\begin{thebibliography}{10}
\expandafter\ifx\csname url\endcsname\relax
  \def\url#1{\texttt{#1}}\fi
\expandafter\ifx\csname urlprefix\endcsname\relax\def\urlprefix{URL }\fi
\expandafter\ifx\csname href\endcsname\relax
  \def\href#1#2{#2} \def\path#1{#1}\fi

\bibitem{bonaventura:2018}
L.~Bonaventura, E.~D. Fern\'andez-Nieto, J.~Garres-D\'iaz, G.~Narbona-Reina,
  \href{https://www.sciencedirect.com/science/article/pii/S0021999118301694}{Multilayer
  shallow water models with locally variable number of layers and semi-implicit
  time discretization}, Journal of Computational Physics 364 (2018) 209 -- 234.
\newline\urlprefix\url{https://www.sciencedirect.com/science/article/pii/S0021999118301694}

\bibitem{fernandezNieto:2016}
E.~D. Fern\'andez-Nieto, J.~Garres-D\'iaz, A.~Mangeney, G.~Narbona-Reina, A
  multilayer shallow model for dry granular flows with the $\mu({I})$ rheology:
  Application to granular collapse on erodible beds, Journal of Fluid Mechanics
  798 (2016) 643--681.

\bibitem{SainteMarie}
J.~Sainte-Marie, Vertically averaged models for the free surface
  non-hydrostatic euler system: derivation and kinetic interpretation, Math.
  Models Methods Appl. Sci. 21~(3) (2011) 459--490.

\bibitem{cocshu}
B.~Cockburn, S.~Chi-Wang, The local discontinuous galerkin method for
  time-dependent convection-diffusion systems, SIAM J. Numer. Anal. 35~(6)
  (1998) 2440--2463.

\bibitem{cockburn}
D.~N. Arnold, F.~Brezzi, B.~Cockburn, L.~D. Marini, Unified analysis of
  discontinuous galerkin methods for elliptic problems, SIAM J. Numer. Anal.
  39~(5) (2002) 1749--1779.

\bibitem{pietro}
D.~A. Di~Pietro, A.~Ern, Mathematical Aspects of Discontinuous Galerkin
  Methods, Vol.~69 of Math{\'e}matiques $\&$ Applications, Springer-Verlag,
  Berlin, 2011.

\bibitem{botasso}
C.~L. Bottasso, S.~Micheletti, R.~Sacco, The discontinuous petrov-galerkin
  method for elliptic problems, Comput. Methods Appl. Mech. Engrg. 191~(31)
  (2002) 3391--3409.

\bibitem{causin}
P.~Causin, R.~Sacco, A discontinuous petrov-galerkin method with lagrangian
  multipliers for second order elliptic problems, SIAM J. Numer. Anal. 43~(1)
  (2005) 280--302.

\bibitem{demko1}
L.~Demkowicz, J.~Gopalakrishnan, A class of discontinuous petrov-galerkin
  methods. part i: The transport equation, Comput. Methods Appl. Mech. Engrg.
  199~(23-24) (2010) 1558--1572.

\bibitem{demko2}
L.~Demkowicz, J.~Gopalakrishnan, Analysis of the dpg method for the poisson
  problem, SIAM J. Numer. Anal. 49~(5) (2011) 1788--1809.

\bibitem{demko3}
L.~Demkowicz, J.~Gopalakrishnan, A class of discontinuous petrov-galerkin
  methods. part ii: Optimal test functions, Numer. Meth. Part. Diff. Equa.
  27~(1) (2011) 70--105.

\bibitem{demko}
L.~Demkowicz, J.~Gopalakrishnan, Discontinuous petrov-galerkin (dpg) method,
  ICES Report 15-20 (2015) 1--21.

\bibitem{brezis}
H.~Brezis, Functional Analysis, Sobolev Spaces and Partial Differential
  Equations, Springer-Verlag, New York, 2011.

\bibitem{ernguermond}
A.~Ern, J.-L. Guermond, Theory and practice of finite element methods, Vol. 159
  of Applied Mathematical Sciences, Springer-Verlag, New York, 2004.

\bibitem{bernardi}
C.~Bernardi, Y.~Maday, F.~Rapetti, Discr{\'e}tisations variationnelles de
  probl{\`e}mes aux limites elliptiques, Vol.~45 of Math{\'e}matiques $\&$
  Applications, Springer-Verlag, Berlin, 2004.

\bibitem{clement}
P.~Cl\'ement, Approximation by finite element functions using local
  regularization, Rev. Fran{\c c}aise Automat. Informat. Recherche
  Op\'erationnelle S\'er. Rouge Anal. Num\'er. 9~(R-2) (1975) 77--84.

\bibitem{freefem++}
F.~Hecht, S.~Auliac, O.~Pironneau, J.~Morice, A.~Le~Hyaric, K.~Ohtsuka,
  P.~Jolivet,
  \href{http://www.freefem.org/ff++/ftp/freefem++doc.pdf}{Freefem++, third
  edition, version 3.58-1} (2018).
\newline\urlprefix\url{http://www.freefem.org/ff++/ftp/freefem++doc.pdf}

\end{thebibliography}

\end{document}